\documentclass[12pt]{amsart}
\usepackage[margin=1in]{geometry}
\usepackage[tt=false]{libertine}
\usepackage{mathabx}
\usepackage[T1]{fontenc}
\usepackage[natbib]{biblatex}

\usepackage[normalem]{ulem}
\usepackage[dvipsnames]{xcolor}
\usepackage[
  colorlinks=true,
  linkcolor=magenta,
  citecolor=magenta,
  urlcolor=magenta
]{hyperref}
\usepackage{enumerate}

\usepackage{tikz}
\usetikzlibrary{calc} 
\usepackage{pgfplots}
\pgfplotsset{compat=1.18}

\newcommand{\R}{\mathbb{R}}
\newcommand{\de}{\delta}
\newcommand{\diam}{\mathrm{diam}}

\newtheorem{theorem}{Theorem}
\newtheorem{proposition}[theorem]{Proposition}
\newtheorem{question}[theorem]{Question}
\newtheorem{lemma}[theorem]{Lemma}
\newtheorem{corollary}[theorem]{Corollary}
\newtheorem{conjecture}[theorem]{Conjecture}

\theoremstyle{definition}

\newtheorem{exercise}[theorem]{Exercise}
\newtheorem{definition}[theorem]{Definition}
\newtheorem{remark}[theorem]{Remark}
\newtheorem{example}[theorem]{Example}

\renewcommand{\hat}{\widehat}
\renewcommand{\tilde}{\widetilde}

\title[Pramanik Lecture Notes]{
Patterns, Distances, and Structure in Sparse Sets}
\author{Paige Bright}

\begin{document}

\begin{abstract}
A central theme in modern harmonic analysis, with strong connections to Ramsey theory, is the tension between randomness and structure: to what
extent can a "large" set avoid prescribed configurations, and when must
such patterns inevitably appear? The answer depends delicately on how size
is measured, with notions such as Hausdorff dimension, Fourier dimension,
and density leading to strikingly different phenomena.

In this mini-course, we explore this theme through two complementary
perspectives. On one hand, we study translation-invariant linear patterns,
a natural class of Euclidean Ramsey-type configurations. We show that sets
of full Hausdorff dimension, and even sets of large Fourier dimension, can
be constructed to avoid many such patterns, highlighting a surprising
degree of flexibility. However, this flexibility breaks down under
Fourier-analytic constraints: sets of sufficiently large Fourier dimension
are forced to contain nontrivial linear patterns, revealing a rigidity
phenomenon tied to Fourier decay.

On the other hand, we investigate distance sets of sparse subsets of
Euclidean space, which encode metric configurations of pairs of points.
Classical results show that sufficiently large sets generate rich families
of distances. We present recent progress demonstrating that even
Lebesgue-null sets with large Hausdorff content must contain intervals of
distances, along with results describing the structure and distribution of
these intervals, and conditions under which sparse unbounded sets realize
all sufficiently large distances.

Together, these results illustrate a unifying principle: while sparse sets
may avoid individual configurations, sufficiently strong analytic or
geometric conditions enforce the emergence of Ramsey-type structure. The
course will emphasize key ideas, methods, and examples, aiming to provide
an accessible entry point to current developments at the interface of
harmonic analysis, geometric measure theory, and Euclidean Ramsey theory. Lectures are by Malabika Pramanik.
\end{abstract}

\maketitle

\section{Introduction}

\begin{question} \label{question1}
When does a set contain a prescribed configuration? In particular, when does largeness force structure to exist? Conversely, how large can a set be while still avoiding a specified pattern?
\end{question}

\noindent The primary aim of these lectures is to provide some partial answers to these questions, both in the positive and negative. This problem is a central one in fractal geometry, but it'd be remiss to not point out that this is a recurring theme in mathematics: sets that are large cannot be completely arbitrary.

\begin{itemize}
    \item If a graph is large enough, certain subgraphs must exist within it.
    \item If $E \subset \mathbb{N}$ is dense, $E$ contains arbitrarily long arithmetic progressions (see Theorem \ref{szemeredi}). 
    \item If a measurable set $E\subset \R^n$ has positive Lebesgue measure, $E$ contains an affine copy of every finite point configuration (see Theorem \ref{steinhaus}).
\end{itemize}

\noindent Let us now carefully define what we mean by a configuration and size.

\begin{definition}
A \textit{configuration}, or \textit{pattern}, is a prescribed geometric or algebraic set.
\end{definition}

For instance, one can consider the following patterns:
\begin{itemize}
    \item The vertices of an equilateral triangle or a square.
    \item A 3-term arithmetic progression (\textit{AP}) of the form $\{x-r, x, x+r\}$ for some $r\neq 0$.
    \item Algebraic solutions to an equation, e.g. $\{x_1,\dots, x_k\}$ such that $P(x_1,\dots, x_k) = 0$ for some polynomial $P$.
\end{itemize}
In this lecture series, we will be exploring when large sets \textit{must have} or \textit{can avoid} such patterns. But what do we mean by a ``large'' set? There are many possible notions of size one could consider, such as positive Lebesgue measure, positive upper density, dimension, and many more. For the most part we will focus our attention to different notions of dimension, but before let us start with a motivating example over the natural numbers.

\begin{theorem}[Szemer\'edi's Theorem \cite{Szemerdi1975}] \label{szemeredi}
    If $E\subseteq \mathbb{N}$ has positive upper density, i.e.
    \[
    \overline{\delta}(E) := \limsup_{N\to \infty} \frac{\#(E \cap \{1,\dots, N\})}{N} >0,
    \]
    then for all $k\geq 3$, $E$ contains a $k$-term arithmetic progression. That is, $E$ contains a collection of points of the form $\{a, a+ d, a + 2d,\dots, a + (k-1) d\}$ for some $d\neq 0$.
\end{theorem}

Using the language we were using before, this shows that if $E$ is \textit{large} in the sense of having positive upper density, it must have some \textit{structure} in the sense of having arbitrarily long arithmetic progressions. Equally important to Szemer\'edi's theorem are results in the opposite direction. Take, for instance, the following result of Salem--Spencer and Behrend:

\begin{theorem} [\cite{Salem1942,Behrend1946}] \label{Behrend}
    For all $N$ sufficiently large, there exists $B_N \subseteq \{0,1,\dots, N-1\}$ with $\#B_N \geq N^{1- o(1)}$ such that $B_N$ avoids 3-term arithmetic progressions.
\end{theorem}

Our interest is to explore analogues of these results in Euclidean space (some of which utilize such discrete statements, see Theorem \ref{ctmBehrend}). To this end, we must first state Steinhaus' theorem, whose proof was covered in Kova\v{c}'s lectures at the same summer school (see \href{https://web.math.pmf.unizg.hr/~vjekovac/files/Lecture_notes_Simons_Budapest.pdf}{here}).

\begin{theorem}[Steinhaus' Theorem \cite{Steinhaus1920}] \label{steinhaus}
Every set of positive Lebesgue measure in $\R^n$ contains an affine copy of every finite set. That is, given $E\subset \R^n$ with $|E|>0$ and any finite set of points $P \subset \R^n$, there exists an affine copy of $P$ within $E.$
\end{theorem}

This is a \textit{rigidity} statement: \textit{every} positive measure set must \textit{necessarily} contain any finite point configuration (up to an affine transformation). Moreover, given Steinhaus' theorem, we necessarily turn our attention to structure within measure zero sets. Still, we will need a corresponding notion of size. As such, we now define \textit{Hausdorff dimension}.

\subsection{Hausdorff Dimension}

Hausdorff dimension will be defined using a family of measures that serve as a fractal analogue (and in fact, generalization) of Lebesgue measure.

\begin{definition}[Hausdorff measure and dimension]
    Let $(X,d)$ be a separable metric space, $E\subset X$, and $\de>0$. Then, we define \begin{itemize}
    \item the \textit{$s$-dimensional Hausdorff content at scale $\de$} via
    \[
    \mathcal H^s_\de(E) = \inf \left\{\sum_i r_i^s : E\subset \bigcup B_i \text{ such that } r_i:= \mathrm{diam}(B_i) < \de\right\},
    \]
    \item the \textit{$s$-dimensional Hausdorff measure} of $E$ as  $\mathcal H^s(E) = \lim\limits_{\de \searrow 0} \mathcal H^s_\de(E)$, and
    \item the \textit{Hausdorff dimension} of $E$ as 
    \[
    \dim_{\mathbb{H}} E = \inf\{s : \mathcal H^s(E) = 0\} = \sup\{s : \mathcal H^s(E) = \infty\}.
    \]
    \end{itemize}
\end{definition}

\begin{remark}
    Note that separability is used to ensure that an infimal \textit{countable} covering can be obtained in the definition of $\mathcal H^s_\de$.
\end{remark}

One can directly check that when $X = \R^n$, the $\mathcal H^n$-measure of $E$ agrees with its Lebesgue measure. Strictly speaking, this is true up to appropriate dimensional constants that should be put into the definition of $\mathcal H^s_\de$ so that the unit cube has volume 1. However, this is unimportant from the perspective of defining Hausdorff dimension. 

\begin{exercise} \label{dimensionbasics}
    Show that if $E\subset \R^n$ is measurable and $\mathcal H^n(E)>0$, then $\dim_{\mathbb{H}} E = n$. Moreover, convince yourself that given a set $E\subset \R^n$ with $\dim E < n$, $\mathcal H^n(E) = 0.$
\end{exercise}

\begin{exercise}
    Show that $\dim_{\mathbb{H}} (\mathbb{Q}) = 0$. Moreover, generalize this to any countable set.
\end{exercise}

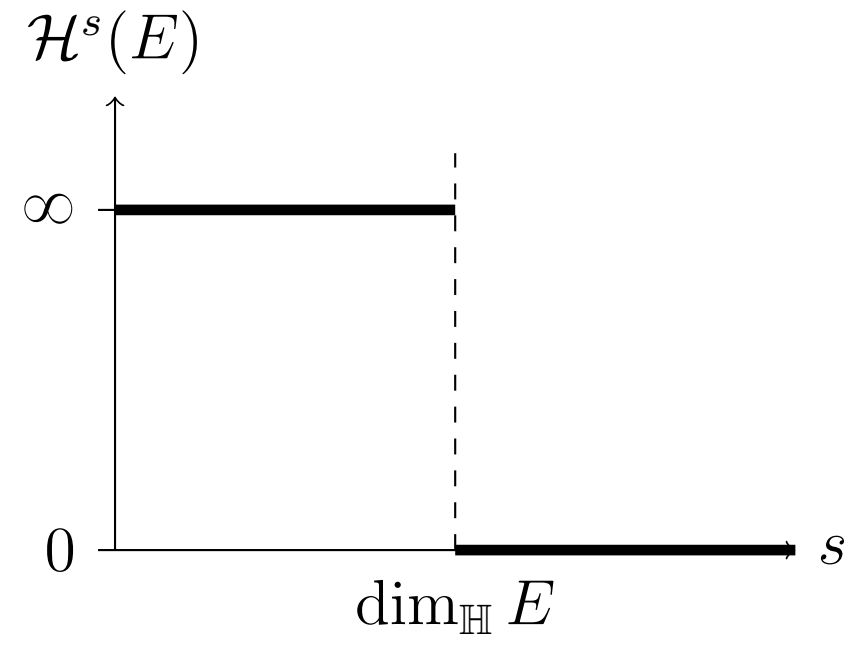
\begin{figure}[h] 
    \centering
    \begin{tikzpicture}[scale=.75]

% axes
\draw[->] (0,0) -- (6,0) node[right] {$s$};
\draw[->] (0,0) -- (0,4) node[above] {$\mathcal{H}^s(E)$};

% dimension point
\def\dimE{3}

% dashed vertical line at dim A
\draw[dashed] (\dimE,0) -- (\dimE,3.5);
\node at (\dimE,-0.5) {$\dim_{\mathbb{H}} E$};

% y-axis ticks
\draw (0,0) -- (-0.15,0);
\node[left] at (-0.15,0) {$0$};

\draw (0,3) -- (-0.15,3);
\node[left] at (-0.15,3) {$\infty$};

% step behavior (visual guide)
\draw[line width=2pt] (0,3) -- (\dimE,3);
\draw[line width=2pt] (\dimE,0) -- (6,0);

\end{tikzpicture}

    \caption{A depiction of $\mathcal H^s(E)$, whose transition point is $\dim_{\mathbb{H}} E$.}
    \label{figure1}
\end{figure}

 To see that Hausdorff dimension is well-defined, note the following proposition (depicted pictorially in Figure \ref{figure1}).

\begin{proposition} \label{prop1}
    Let $(X,d)$ be a metric space, $E\subset X$, and $0 \leq s < t < \infty$. Then, 
    \begin{itemize}
        \item If $\mathcal H^s(E) < \infty$, then $\mathcal H^t (E) = 0$.
        \item If $\mathcal H^t(E) >0$, then $\mathcal H^s(E) = \infty.$
    \end{itemize}
\end{proposition}

\begin{proof}
Let $\de >0$ and let $E\subset \bigcup B_i$ with $\diam(B_i) \leq \de$. Then, by definition, 
\[
\sum_i \diam(B_i)^t \leq \sum_i \de^{t-s}\diam(B_i)^s.
\]
Taking infimums, it follows that 
\[
\mathcal H^t_\de(E) \leq \de^{t-s} \mathcal H_\de^s(E).
\]
The result follows immediately from noting $\de^{t-s} \to 0$ as $\de \searrow 0$.
\end{proof}

Hausdorff dimension agrees with other reasonable notions of topological dimension. For instance, the Hausdorff dimension of a curve is 1, and that of a surface is 2. However, Hausdorff dimension can also take non-integer values.

\begin{example} \label{ex:midthird}
    The classical middle thirds Cantor set, $\mathcal C_{1/3}$, has Lebesgue measure zero and Hausdorff dimension $\frac{\log 2}{\log 3}$.
\end{example}

\begin{exercise}
    Sketch a proof of Example \ref{ex:midthird} by considering what the minimal covering in $\mathcal H^s_{3^{-k}}(\mathcal C_{1/3})$ should be for $k\geq 1$ and arbitrary $s$.
\end{exercise}

\begin{exercise} \label{logjlogk}
    Based on the above example, use a Cantor-like construction to construct sets of Hausdorff dimension $\frac{\log j}{\log k}$ for $j,k\geq 2$.
\end{exercise}

\begin{exercise}
    Show that for any countable collection of Borel sets $E_i \subset \R^n$, 
    \[
    \dim_{\mathbb{H}} \bigcup E_i = \sup_i \dim_{\mathbb{H}} E_i.
    \]
    Using that the Cantor sets in Exercise \ref{logjlogk} have Lebesgue measure 0 by Exercise \ref{dimensionbasics}, construct a set $E\subset \R^1$ with null Lebesgue measure such that $\dim_{\mathbb{H}} E = 1$. Generalize this to $\R^n.$
\end{exercise}

The above definition of Hausdorff dimension can be quite difficult to use in practice; requiring finding families of infimal coverings of a set. To this end, we note Frostman's lemma, which will be of more practical use throughout these lectures.

\begin{theorem}[Frostman's lemma, cf. {\cite[Theorem 8.8]{Mattila95}}] \label{Frostman}
    Given $E\subset \R^n$ Borel, we have 
    \[
    \dim_{\mathbb{H}} (E) = \sup\{s \in [0,n] : \exists \text{ a probability measure } \mu, \mathrm{supp}(\mu) \subset E, \sup_{x\in \R^n,\, r>0} \frac{\mu(B(x,r))}{r^s} < \infty\}.
    \]
    Such a measure is called a \textit{Frostman measure}.
\end{theorem}

One may reasonably ask what such a measure $\mu$ has to do with the geometry of our set $E$. It is worth thinking about what this (heuristically) implies about the number of $\de$-balls needed to cover $E$ for sufficiently small $\de$. It may help to think about this when $E$ is a singleton, and thus $\mu$ must be the Dirac delta measure.

We are now able to revise our original Question \ref{question1}.

\begin{question} \label{revisedquestion}
    Must a set of large Hausdorff dimension contain a prescribed configuration?
\end{question}

To this end, let us state two results of this flavor.

\begin{theorem}[Keleti \cite{Keleti1998,keleti2008}] \label{keleti}
    Given three points $P = \{x,y,z\}$ (in fact, any finite set of points), there exists $E\subset [0,1]$ of full Hausdorff dimension (i.e. $\dim_{\mathbb{H}} E = 1$) avoiding all affine copies of $P$.
\end{theorem}

\begin{theorem}[M\'ath\'e \cite{Mth2011CoveringTR}]
    There exists a compact set $E \subset [0,1]$ with $\dim_{\mathbb{H}}(E) = 0$ containing an affine copy of every finite set.
\end{theorem}

\begin{remark}
    Note that these results imply that the answer to Question \ref{revisedquestion} is, in general, \textit{no}.
\end{remark}

The above two results are \textit{flexibility} statements: A set can be large or small and contain/not contain a prescribed configuration. Throughout these lectures, we will try to highlight the tension between \textit{rigidity} and \textit{flexibility} results.

\subsection{A Set Avoiding 3-APs}

Motivated by Theorem \ref{keleti}, we construct a set of full Hausdorff dimension avoiding 3-term APs. In Section \ref{subsec:Behrend}, we work on many more such constructions.

\begin{theorem} \label{ctmBehrend}
    There exists a set $E\subset [0,1]$ of full Hausdorff dimension such that $E$ avoids all 3-term arithmetic progressions.
\end{theorem}

\noindent We will demonstrate Theorem \ref{ctmBehrend} through two telling but different methodologies.

\textbf{Method 1: Cantor Constructions}: The \textit{Cantor philosophy} behind \textbf{Method 1} is to construct a nested sequence of sets where some intervals are retained or removed at each scale. We will do so in a way that gives large Hausdorff dimension (which will follow from which intervals are retained), while still avoiding this pattern (which will follow from which intervals are removed).

As stated previously, our set will be constructed iteratively. Let $E_0 = [0,1]$. We will break $E_0$ into $10 N_1$ intervals of length $\frac{1}{10N_1}$ avoiding 3-term arithmetic progressions at that scale. In particular, given $N_1$ sufficiently large and $B_{N_1}$ as in Theorem \ref{Behrend}, let 
\[
E_1 = \{\bigcup_{j_1} \frac{10 j_1}{10 N_1} + [0,\frac{1}{10N_1}] : j_1 \in B_{N_1}\}.
\]
Note that while $E_1$ clearly contains 3-APs (as it contains intervals), by our construction any such AP must lie in exactly one of the $\frac{1}{10 N_1}$ intervals in $E_1.$ More generally, we will let $E_k$ be obtained by decomposing every basic interval of $E_{k-1}$ into $10 N_k$ intervals of length $\frac{1}{10N_k}$ equal pieces, and putting a copy of the Behrend set ($10 B_{N_k}$) within it. Finally, we let $E = \bigcap E_k$.

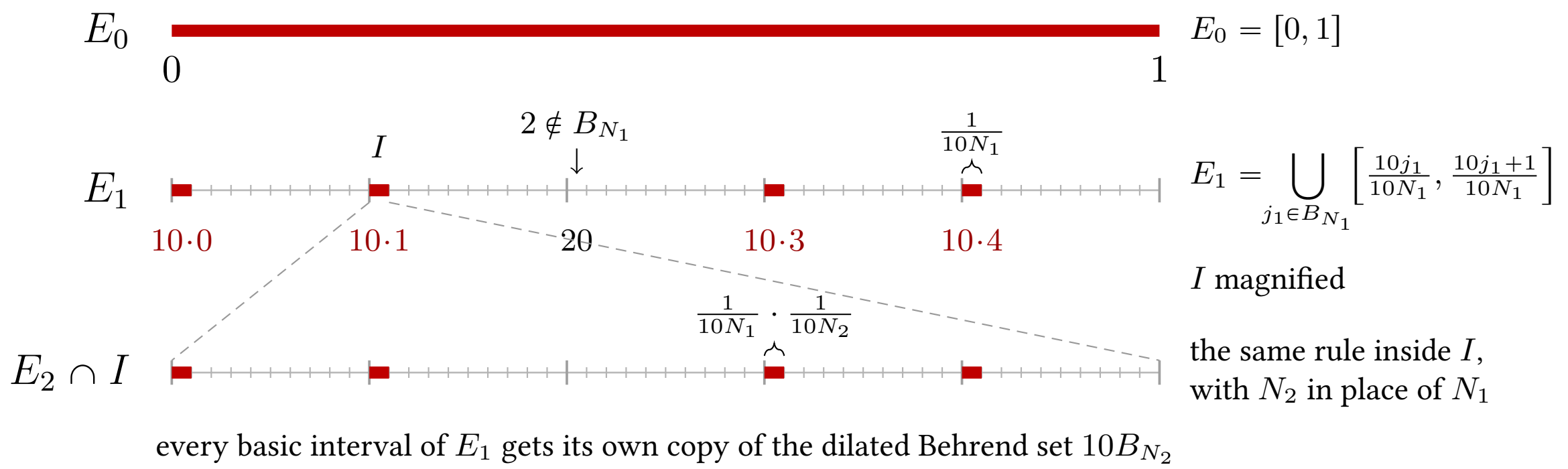
\begin{figure}[h]
    \centering
    \begin{tikzpicture}[scale=.75,
    axis/.style={line width=0.5pt, gray!55},
    kept/.style={line width=3.4pt, red!75!black},
    tick/.style={gray!60, line width=0.4pt},
    big/.style={gray!75, line width=0.7pt},
    note/.style={font=\scriptsize, black},
  ]
  % ---- toy parameters: N = 5, Behrend set B = {0,1,3,4} (3-AP free) ----
  \def\L{13}
  \def\N{5}
  \def\B{0,1,3,4}
  \pgfmathsetmacro{\blk}{\L/\N}          % block of ten subintervals
  \pgfmathsetmacro{\sub}{\L/(10*\N)}     % one subinterval = 1/(10 N_1)
  \def\Iz{1}                             % which retained interval gets zoomed
 
  % ===================== E_0 =====================
  \node[anchor=east] at (-0.35,0) {$E_0$};
  \draw[kept] (0,0) -- (\L,0);
  \node[below=2pt] at (0,-0.05) {\small $0$};
  \node[below=2pt] at (\L,-0.05) {\small $1$};
  \node[note, anchor=west] at (\L+0.2,0) {$E_0=[0,1]$};
 
  % ===================== E_1 =====================
  \begin{scope}[yshift=-2.1cm]
    \node[anchor=east] at (-0.35,0) {$E_1$};
    \draw[axis] (0,0) -- (\L,0);
    \foreach \j in {0,...,50} { \draw[tick] ({\j*\sub},-0.07) -- ({\j*\sub},0.07); }
    \foreach \b in {0,...,\N} { \draw[big] ({\b*\blk},-0.16) -- ({\b*\blk},0.16); }
    % retained subintervals: index 10 j_1, j_1 in B'_{N_1}
    \foreach \b in \B {
      \draw[kept] ({\b*\blk},0) -- ({\b*\blk+\sub},0);
      \node[below=7pt, font=\scriptsize, red!60!black] at ({\b*\blk+0.5*\sub},0) {$10\!\cdot\!\b$};
    }
    % the removed block
    \node[below=7pt, font=\scriptsize, black] at ({2*\blk+0.5*\sub},0) {$20$};
    \node[note, align=center] at ({2*\blk+0.5*\sub},0.85) {$2\notin B_{N_1}$};
    \draw[black, line width=0.4pt, ->, shorten >=2pt]
      ({2*\blk+0.5*\sub},0.55) -- ({2*\blk+0.5*\sub},0.12);
    % length of a basic interval
    \draw[decorate, decoration={brace, amplitude=3.5pt}]
      ({4*\blk},0.22) -- ({4*\blk+\sub},0.22)
      node[midway, above=1pt, font=\scriptsize] {$\tfrac{1}{10N_1}$};
    \node[note, anchor=west, align=left] at (\L+0.2,0)
      {$\displaystyle E_1=\!\!\bigcup_{j_1\in B_{N_1}}\!\!
        \Big[\tfrac{10j_1}{10N_1},\tfrac{10j_1+1}{10N_1}\Big]$};
    % name the interval that gets magnified, above its red bar
    \node[font=\scriptsize, above=5pt] at ({\Iz*\blk+0.5*\sub},0) {$I$};
  \end{scope}
 
  % ===================== zoom =====================
  \draw[densely dashed, black!40] ({\Iz*\blk},-2.26) -- (0,-4.34);
  \draw[densely dashed, black!40] ({\Iz*\blk+\sub},-2.26) -- (\L,-4.34);
  \node[note, anchor=west] at (\L+0.2,-3.3) {$I$ magnified};
 
  % ===================== E_2 inside I =====================
  \begin{scope}[yshift=-4.5cm]
    \node[anchor=east, align=right] at (-0.35,0) {$E_2\cap I$};
    \draw[axis] (0,0) -- (\L,0);
    \foreach \j in {0,...,50} { \draw[tick] ({\j*\sub},-0.07) -- ({\j*\sub},0.07); }
    \foreach \b in {0,...,\N} { \draw[big] ({\b*\blk},-0.16) -- ({\b*\blk},0.16); }
    \foreach \b in \B { \draw[kept] ({\b*\blk},0) -- ({\b*\blk+\sub},0); }
    \draw[decorate, decoration={brace, amplitude=3.5pt}]
      ({3*\blk},0.22) -- ({3*\blk+\sub},0.22)
      node[midway, above=1pt, font=\scriptsize] {$\tfrac{1}{10N_1}\cdot\tfrac{1}{10N_2}$};
    \node[note, anchor=west, align=left] at (\L+0.2,0)
      {the same rule inside $I$,\\ with $N_2$ in place of $N_1$};
    \node[note, below=13pt] at (\L/2,0)
      {every basic interval of $E_1$ gets its own copy of the dilated Behrend set $10B_{N_2}$};
  \end{scope}
\end{tikzpicture}
    \caption{The construction in \textbf{Method 1} with $N_k=5$ and $B_{N_k}=\{0,1,3,4\}$,
     which is $3$-AP free. Iterating gives $E=\bigcap E_k$.}
    \label{fig:Method1}
\end{figure}

\begin{remark}
    Here, one can choose 
    \[
    100 \leq N_1 \leq N_2 \leq \dots
    \]
    growing sufficiently slowly. However, for other similar constructions we will need to be more careful about this choice of $N_k$.
\end{remark}

\begin{exercise}
    Given $\#B_N \geq N^{1- o(1)}$, show $\dim_{\mathbb{H}} E = 1.$
    Moreover, show that $E$ is 3-AP free.
\end{exercise}

There are some immediate pros and cons to \textbf{Method 1}. A big pro is that we can easily convert a discrete example from additive combinatorics into a fractal one with the appropriate dimension. On the other hand, a big con is that we need have said discrete example in hand. Keleti was able to bypass these particular pros and cons by the following methodology. 

\textbf{Method 2: Cantor Constructions with Memory}: Cantor constructions hardly ever consider what happened at the previous stage(s) in the construction to inform the current stage. In this method, we present a construction of a Cantor set with ``memory'' due to Keleti that does precisely this. We begin as before with $E_0 = [0,1]$, and let 
\[
E_1 = \bigcup_{j_1} \left\{\frac{j_1}{6N_1} + [0,\frac{1}{6N_1}] : j_1 = 0 \pmod 6\right\}.
\]
Call the basic intervals in the above union, listed in order, $J_1^{(1)}, J_2^{(1)},\dots$, and let 
\[
Q_1 = \{J_1^{(1)}, J_2^{(1)},\dots, \}
\]
be an ordered list of all basic intervals in $E_1.$ We now define $E_2$ as 
\[
E_2 = \bigcup_{j_1,j_2} \{I_{j_1,j_2} := \frac{j_1}{6N_1} + \frac{j_2}{6^2N_1N_2} + [0, \frac{1}{6^2 N_1 N_2}]\}
\]
where if $I_{j_1,j_2} \subseteq J_1^{(1)}$, i.e. the first element of the queue in $Q_1$, we pick $j_2 = 3 \pmod 6$, and otherwise pick $j_2 = 0 \pmod 6$. We then update $Q_1$ so that ``no one loses their place in the queue," i.e. 
\[
Q_2 = \{J_1^{(1)},J_2^{(1)},\dots, J_1^{(2)},J_2^{(2)},\dots\}
\]
where $J_i^{(2)}$ are the basic intervals in $E_2$. In general, to construct $E_k$, we decompose every basic interval $I$ in $E_{k-1}$ into $6N_k$ equal parts. If $I$ is contained in the $k$-th member of the queue $Q_{k-1}$, we choose the subinterval indexed by $j_k = 3 \pmod 6$, and $j_k = 0 \pmod 6$ otherwise. Finally, we let 
\[
Q_k = Q_{k-1}
 \cup \{\text{the basic intervals in $E_k$ (ordered)}\}.
 \]

We claim that $E = \bigcap E_k$ is 3-AP free. Indeed, suppose $x, y, z \in E$ with $x < y < z$ and $x + z = 2y$. Since $x \neq y$, there is a basic interval $J$ (at some stage) containing $x$ but not $y$; as $J$ is an interval and $z > y$, it does not contain $z$ either. Since the queue is only ever appended to, $J$ is the $k$-th member of $Q_{k-1}$ for some $k$, and $J$ is a basic interval of some stage strictly before $k$. Consequently, at stage $k$, the basic interval of $E_{k-1}$ containing $x$ lies inside $J$, while those containing $y$ and $z$ do not.

Now let $h := (6^k N_1 \cdots N_k)^{-1}$ denote the length of a basic interval of $E_k$. For each $w \in \{x, y, z\}$, we may write
\[
w = a_w + j_w h + \epsilon_w,
\]
where $a_w$ is the left endpoint of the basic interval of $E_{k-1}$ containing $w$, $j_w \in \{0, \ldots, 6N_k - 1\}$ is the index chosen at stage $k$, and $\epsilon_w \in [0, h)$. Note that $\epsilon_w < h$ strictly, since the last subinterval (with index $6N_{k+1} - 1 \equiv 5 \pmod 6$) is never kept at stage $k+1$. By the above, $j_x \equiv 3 \pmod 6$ while $j_y \equiv j_z \equiv 0 \pmod 6$. Moreover, every left endpoint $a_w$ is an integer multiple of $6N_k h$, the length of a basic interval of $E_{k-1}$. Substituting into $x + z - 2y = 0$ and dividing by $h$, we find
\[
\frac{a_x + a_z - 2a_y}{h} + (j_x + j_z - 2j_y) = -\frac{\epsilon_x + \epsilon_z - 2\epsilon_y}{h} \in (-2, 2).
\]
The left-hand side is an integer congruent to $j_x + j_z - 2j_y \equiv 3 \pmod 6$, since the first term is a multiple of $6N_k$. However, the only integers in $(-2, 2)$ are $-1$, $0$, and $1$, none of which is congruent to $3$ modulo $6$. This is a contradiction, so $E$ contains no 3-APs.

\begin{exercise} \label{exercise:h-dimBehrend}
    Choosing $N_1,\dots, N_k$ such that 
    \[
    \lim_{k\to \infty} \frac{\log (N_1\cdots N_k)}{ \log(6^k N_1\cdots N_k)} = 1,
    \]
    show that $\dim_{\mathbb{H}} E = 1$.
\end{exercise}

As before, we note some pros and cons to \textbf{Method 2}. Firstly, a big pro is that it avoids the need to consider an analogous construction over the integers as we saw in \textbf{Method 1}. On the other hand, it can be hard to adapt the construction in \textbf{Method 2} to sets that have large \textit{Fourier dimension}. This notion of size will be the focus of Section \ref{lec2}.

\subsection{Behrend-type Constructions} \label{subsec:Behrend}

We now discuss generalizations of \textbf{Method 1}, i.e. Behrend-type constructions, following the work of Liang--Pramanik \cite{LiangPramanik}.

The main theorem we will discuss involves sets that avoid nontrivial solutions to translation invariant linear equations (sometimes referred to as TILEs). These are equations in the set
\[
\mathcal F := \left\{f(x_0, \dots, x_v) = m_0 x_0 - \sum_{i=1}^v m_i x_i : m_i \in \mathbb{N}, m_0 = \sum_{i=1}^v m_i, v\geq 2\right\}.
\]
Note that additionally, without loss of generality, we may assume 
\[
\mathrm{gcd}(m_0,m_1,\dots, m_v) = 1
\]
as otherwise we can factor out the greatest common divisor in the above equations.

Before stating the main result, let us note a few things. Firstly, we call these equations translation invariant as solutions are invariant under affine maps. That is, given any affine transformation $T: \R\to \R$, 
\[
f(x_0,\dots, x_v) = 0 \quad \text{ if and only if } \quad f(Tx_0,\dots, Tx_v) = 0.
\]
Moreover, by the translation invariance of functions $f\in \mathcal F$, it follows that for any $x\in \R$, $f(x,x,\dots, x) = 0$. Thus, we will say that a solution $x$ is \textit{nontrivial} if the entries $x_i$ are not all equal. Finally, we note that 3-term arithmetic progressions are nontrivial solutions to the function $f\in \mathcal F$ of the form 
\[
f(x_0,x_1, x_2) = 2x_0 - (x_1 + x_2).
\]
In this sense, the following theorem generalizes Theorem \ref{ctmBehrend}.

\begin{theorem}[Liang--Pramanik {\cite[Theorem 1.5]{LiangPramanik}}] \label{liangpramanik}
    Given any finite $\mathcal G \subset \mathcal F$, there exists a set $E\subset [0,1]$ of full Hausdorff dimension such that $E$ avoids all nontrivial solutions to functions in $\mathcal G$.
\end{theorem}

\begin{remark}
    In fact, one can prove the existence of such a set of full \textit{Fourier} dimension, see Section \ref{sec:randomness}.
\end{remark}

Our proof of Theorem \ref{liangpramanik} will follow the Behrend-like construction we saw previously. In particular, we begin by proving the following discrete statement.

\begin{proposition}[Liang--Pramanik {\cite[Proposition 4.2]{LiangPramanik}}] \label{integeravoidance}
    Fix $\mathcal G \subseteq \mathcal F$ finite. Let 
    \[
    A > 2\max\{m_0 : \exists f\in \mathcal G \text{ with } f(x_0,\dots, x_v) = m_0 x_0 - \sum_{i=1}^v m_i x_i\}.
    \]
    Then, for all sufficiently large $N\geq N_0(A)$, there exists a set $Y_N = Y_N(\mathcal G) \subset \{1,\dots, N\}$ with 
    \[
    \# Y_N \geq N^{1-o(1)}
    \]
    avoiding all nontrivial solutions to functions in $\mathcal G$.
\end{proposition}

\begin{proof}
    Fix large numbers $n,d$ to be determined shortly. For every $0\leq y \leq (Ad)^n - 1$ we use a base $Ad$ expansion to write 
    \[
    y = y_0 + y_1 (Ad) + y_2 (Ad)^2 + \cdots + y_{n-1}(Ad)^{n-1}
    \]
    for a unique choice of $(y_0,\dots, y_{n-1}) \in \{0,1,\dots, d-1\}^n$. Let 
    \[
    Y[y] = (y_0,\dots, y_{n-1}).
    \]
    Moreover, for all $k\geq 0$, let 
    \[
    S_k(d,n) = \left\{y \in \{0,\dots, (Ad)^{n}\} : \lVert Y[y]\rVert^2 := \sum_{i=0}^{n-1} y_i^2 = k\right\}.
    \]
    We make two claims that will conclude the proof.

    \textbf{Claim 1}: There exists $k_0\in \{0,1,2,\dots, n(d-1)^2\}$ such that $\#S_k(d,n) > d^{n-2}/n$.

    \textbf{Claim 2}: Each $S_k(d,n)$, if nonempty, avoids all nontrivial zeroes of functions in $\mathcal G$.

    Let us first show that these two claims finish the proof. For sufficiently large $N$, let $n = \lfloor \ln N\rfloor$ and let $d\in \mathbb{N}$ such that $(Ad)^n \leq N < (A(d+1))^n$ and then set $Y_N := S_{k_0}(d,n).$ This $Y_N$ will be the desired one in the statement of the theorem. By \textbf{Claim 1}, $S_{k_0}(d,n)$ is nonempty, and thus by \textbf{Claim 2} the set $Y_N$ contains no nontrivial zeroes of any $f\in \mathcal G$ as desired. Moreover, by \textbf{Claim 1}, we see
    \[
    \#Y_N > \frac{d^{n-2}}{n} \geq N^{1-o(1)}.
    \]
    Note that here we use that $N$ is sufficiently large dependent on $A$. This completes the proof.

    We now just need to verify \textbf{Claims 1} and \textbf{2}. As for \textbf{Claim 1}, we proceed by a pigeonholing argument. In particular, by definition of $S_k$, the entries of $Y[y]$ are required to lie in $\{0,\dots, d-1\}$. Since $k = \lVert Y[y]\rVert^2$. This forces $0 \leq k \leq n(d-1)^2$, and the number of possible choices of $k$ such that $S_k(d,n) \neq \emptyset$ is thus at most $n(d-1)^2 +1$. Finally, the sets $S_k(d,n)$ are disjoint, and thus 
    \[
    \sum_k \# S_k = d^n.
    \]
    Thus, by the pigeonhole principle, there exists a $k_0 \in \{0,1,2,\dots, n(d-1)^2\}$ such that
    \[
    \#S_{k_0}(d,n) \geq \frac{d^n}{n(d-1)^2 + 1} > \frac{d^{n-2}}{n},
    \]
    finishing the proof of \textbf{Claim 1}.

    Now onto \textbf{Claim 2}, which will follow by a triangle inequality argument. Suppose there exists 
    \[
    (x_0,\dots, x_v) \in [S_k(d,n)]^{v+1}
    \]
    that is a zero of a function $f\in \mathcal G$ of the form $f(x_0,\dots, x_v) := m_0 x_0 - \sum_{i=1}^v m_i x_i$. As before, let $X[x_i] = (x_{i,0},\dots, x_{i,n-1})$ denote the vector of digits in $x_i$ in base $Ad$. Then, by definition, $\lVert X[x_i]\rVert = \sqrt{k}$, and $0 \leq x_{i,j} < d$ for all $(i,j)$. Using this expansion, we have 
    \[
    \sum_{j=0}^{n-1} (m_0 x_{0,j}) (Ad)^j = \sum_{j=0}^{n-1} \left[\sum_{i=1}^v m_i x_{i,j}\right](Ad)^j.
    \]
    Since $A> 2m_0 = 2(m_1+ \dots + m_v)$ by assumption, it follows that for any $j \in \{0,\dots, n-1\},$ $0 \leq m_0 x_{0,j} < m_0 d < Ad$. Additionally, $0 \leq \sum_{i=1}^v m_i x_{i,j} < d \sum_{i=1}^v m_i < Ad$. The uniqueness of the digit expansion implies that the vectors $\{X[x_i] : 0 \leq i \leq v\}$, which lie on the sphere of radius $\sqrt{k}$ by assumption, additionally satisfy the linear relation 
    \[
    m_0 X[x_0] = \sum_{i=1}^v m_i X[x_i].
    \]
    Therefore, by the triangle inequality
    \[
    m_0 \sqrt{k} = \lVert m_0 X[x_0] \rVert \leq \sum_{i=1}^v m_i \lVert X[x_i]\rVert = m_0 \sqrt{k}.
    \]
    This in fact implies that the above inequalities are equalities, and thus $X[x_i]$ are proportional. Since they all lie on the sphere of radius $
    \sqrt{k}$, it follows that $X[x_1] = \cdots = X[x_v]$ and thus $x_0 = x_1 = x_2 = \dots = x_v$. Hence, any solutions $S_k(d,n)$ contains for some $f\in \mathcal G$ are trivial solutions. This concludes the proof of \textbf{Claim 2} and thus the proposition.
\end{proof}

\begin{exercise}
    Work through Proposition \ref{integeravoidance} with $\mathcal G = \{f(x_0,x_1,x_2) = 2x_0 - (x_1 + x_2)\}$ and deduce Theorem \ref{Behrend}.
\end{exercise}

To conclude the existence of a set $E\subset [0,1]$ of full Hausdorff dimension avoiding all nontrivial solutions to functions in $\mathcal G$, we proceed similarly to the case of Theorem \ref{ctmBehrend} using Proposition \ref{integeravoidance}. In particular, our construction will consist of (an intersection of) a union of intervals $\bigcup_{n} I_n = E$ at smaller and smaller scales such that if $(x_0,\dots, x_v) \in E^{v+1}$ is a solution to some function $f\in \mathcal G$, all of the $x_i$ are contained in a single basic interval.

Fix $N \in \mathbb{N}$. For $j \in \{1, \dots, N-1\}$ let
\[
  I_{N,j} = [\tfrac{j}{N}, \tfrac{j+1}{N}],
\]
and for $Y \subseteq \{1,\dots,N-1\}$ write
\[
  I(Y) = \bigcup_{j \in Y} I_{N,j}.
\]
We say that $I(Y)$ \emph{partially avoids} nontrivial solutions to functions in $\mathcal{G}$
if every solution it contains is confined to a single basic interval: whenever
$x_0, \dots, x_v \in I(Y)$ satisfy $f(x_0,\dots,x_v) = 0$ for some $f \in \mathcal{G}$,
there is a single $j \in Y$ with $x_0, \dots, x_v \in I_{N,j}$.
 
\begin{lemma}[Liang--Pramanik {\cite[Lemma 4.3]{LiangPramanik}}]\label{lem:partial-avoid}
Fix $\mathcal{G} \subset \mathcal{F}$ finite and let $A$ be as in Proposition \ref{integeravoidance}.
For all sufficiently large $N$ there exists $U_N \subseteq \{1, \dots, N-1\}$ with
\[
  \#U_N \geq N^{1-o(1)}
\]
such that $I(U_N)$ partially avoids nontrivial solutions to functions in $\mathcal{G}$.
\end{lemma}
 
\begin{proof}
Let $\widetilde{N} = \lfloor \tfrac{N-1}{A} \rfloor$. For $N$ large, $\widetilde{N} \geq N_0(A)$,
so Proposition \ref{integeravoidance} supplies a set
$Y_{\widetilde{N}} \subseteq \{1, \dots, \widetilde{N}\}$ with
$\# Y_{\widetilde{N}} \geq \widetilde{N}^{1-o(1)}$ that avoids all nontrivial solutions to
functions in $\mathcal{G}$. Put
\[
  U_N = A \cdot Y_{\widetilde{N}} = \{ A j : j \in Y_{\widetilde{N}}\}.
\]
Since $A\widetilde{N} \leq N-1$ we have $U_N \subseteq \{1,\dots,N-1\}$, and since $A = A(\mathcal{G})$
is a constant,
\[
  \# U_N = \#Y_{\widetilde{N}} \geq \widetilde{N}^{\,1-o(1)}
  = (\tfrac{N-1}{A})^{1-o(1)} = N^{1-o(1)}.
\]
 
\noindent Suppose $x_0, \dots, x_v \in I(U_N)$ satisfy $f(x_0,\dots,x_v) = 0$ for some $f\in \mathcal G$ with
$$f(x_0,\dots,x_v) = m_0 x_0 - \sum_{i=1}^{v} m_i x_i.$$
Each $x_i$ lies in a basic interval of $I(U_N)$, i.e. there are $j_i \in Y_{\widetilde{N}}$ and
$\delta_i \in [0,1]$ with
\[
  x_i = \frac{A j_i + \delta_i}{N}, \qquad 0 \leq i \leq v.
\]
Substituting into $f(x_0,\dots,x_v) = 0$ and multiplying by $N$,
\[
  A( m_0 j_0 - \sum_{i=1}^{v} m_i j_i )
  = -( m_0 \delta_0 - \sum_{i=1}^{v} m_i \delta_i ).
\]
The left-hand side is $A$ times an integer. For the right-hand side, translation invariance gives
$\sum_{i=1}^v m_i = m_0$, so both $m_0 \delta_0$ and $\sum_{i=1}^v m_i \delta_i$ lie in $[0, m_0]$
and therefore
\[
  | m_0 \delta_0 - \sum_{i=1}^{v} m_i \delta_i | \leq m_0 < A,
\]
by the choice $A > 2\max\{m_0\}$. An integer multiple of $A$ of absolute value strictly less
than $A$ must vanish, hence
\[
  m_0 j_0 - \sum_{i=1}^{v} m_i j_i = 0,
  \qquad\text{that is,}\qquad
  f(j_0, \dots, j_v) = 0 .
\]
This is a solution to $f$ inside $Y_{\widetilde{N}}$, which avoids all nontrivial ones, so
$j_0 = j_1 = \cdots = j_v =: j$. Consequently every $x_i$ lies in the single basic interval
$I_{N, Aj}$, with $Aj \in U_N$, as claimed.
\end{proof}

\begin{exercise}
    Using a Behrend-type construction and Lemma \ref{lem:partial-avoid}, prove Theorem \ref{liangpramanik}. That is: for any $\mathcal G \subset \mathcal F$ finite, there exists a full Hausdorff dimensional set $E\subset [0,1]$ such that $E$ avoids nontrivial solutions to functions in $\mathcal G$.
\end{exercise}

\noindent In the Section \ref{lec2}, we will improve upon the above construction to find a set $E\subset [0,1]$ of full \textit{Fourier} dimension avoiding nontrivial solutions to functions in $\mathcal G$.

\section{Fourier Dimension} \label{lec2}

The focus of this section will be to motivate and discuss \textit{Fourier dimension} as a notion of size good for detecting/avoiding patterns. However, before discussing this further, let us first mention a result of M\'ath\'e's that will be important for the purposes of contrast throughout this section.

\begin{theorem}[M\'ath\'e \cite{Mathe2017}] \label{matheFulldim}
    There exists a full (Hausdorff) dimensional rationally independent set $E\subset \R$. That is, $\dim_{\mathbb{H}}(E) = 1$, and for all $k\geq 2$, if there exists distinct $x_i \in E$ and $c_i \in \mathbb{Q}$ such that
    \[
    c_1 x_1 + c_2 x_2 + \cdots + c_k x_k = 0,
    \]
    then $c_1=c_2=\cdots =c_k=0$.
\end{theorem}

As we will see, due to a result of Liang--Pramanik \cite{LiangPramanik}, such a set must have Fourier dimension 0. To this end, let us now define the Fourier dimension of a set $E\subset \R^n$ in a similar manner to the definition of Hausdorff dimension in Theorem \ref{Frostman}.

\begin{definition}
    Given $E\subset \R^n$, the \textit{Fourier dimension} of $E$ is defined as 
    \[
    \dim_{\mathbb{F}} E = \sup\{s \in [0,n] : \exists \text{ a probability measure } \mu, \mathrm{supp}(\mu) \subset E, |\hat{\mu}(\xi)| \lesssim \langle \xi\rangle^{- s/2}\}.
    \]
    Note that here, 
    \[
    \hat{\mu}(\xi) = \int_{\R^n} e^{-2\pi i x\cdot \xi} \,\mathrm d\mu(x) \quad \text{ and } \quad \langle \xi\rangle = (1 + |\xi|^2)^{1/2}.
    \]
\end{definition}

\begin{exercise} \label{hatmuLp}
    Let $E\subset \R^n$ be a set with $\dim_{\mathbb{F}}(E) > s$. Show that there exists a probability measure $\mu$ with $\mathrm{supp}(\mu) \subset E$ such that $\hat{\mu} \in L^{2n/s}(\R^n).$
\end{exercise}

\begin{exercise}
    Show that $\dim_{\mathbb{F}}(\{0\}) = 0$ and $\dim_{\mathbb{F}}([0,1]) = 1$ by considering the uniform measures on each. Additionally, show that $\dim_{\mathbb{F}}(\mathcal C_{1/3}) = 0$ where $\mathcal C_{1/3}$ is the middle thirds Cantor set. For proof of the latter, see \cite{Lyons2020}. 
\end{exercise}

\noindent Using an energy-theoretic version of Frostman's lemma, one can aptly show that for all $E\subset \R^n$, 
\begin{equation}\label{FourierleqHausdorff}
\dim_{\mathbb{F}}(E) \leq \dim_{\mathbb{H}}(E).
\end{equation}
We omit the proof, but give the statement as an exercise.

\begin{exercise}
    Let $E\subset \R^n$ and $s<t < \dim_{\mathbb{F}}(E)$. Then, by definition, there exists a measure $\mu$ such that $|\hat{\mu}(\xi)|\lesssim \langle \xi\rangle^{-s/2}$. Show that the Riesz $t$-energy of $\mu$, $I_t(\mu)$, satisfies 
    \[
    I_t(\mu) := \int_{\R^n} |\hat{\mu}(\xi)|^2 |\xi|^{t-n} \,\mathrm d\xi < \infty.
    \]
    A standard energy-theoretic version of Frostman's lemma states that 
    \begin{equation} \label{EnergyFrostman}
    \dim_{\mathbb{H}}(E) = \sup\{s \in [0,n]: \exists \text{ a probability measure } \mu, \mathrm{supp}(\mu)\subset E, I_s(\mu) < \infty\}.
    \end{equation}
    Deduce \eqref{FourierleqHausdorff}.
\end{exercise}

From \eqref{FourierleqHausdorff}, it is reasonable to make the following definition.

\begin{definition}
    We say that a set $E\subset \R^n$ is \textit{Salem} if $\dim_{\mathbb{F}}(E) = \dim_{\mathbb{H}}(E)$.
\end{definition}

The notion of Fourier dimension begs the following question akin to that of Question \ref{revisedquestion}.

\begin{question}
    Must a set of large Fourier dimension contain a prescribed configuration?
\end{question}

\noindent As before, let us state two results of this flavor.

\begin{theorem}[Shmerkin {\cite[Theorem B]{Shmerkin2016}}] \label{shmerkin3AP}
    There exists a full dimensional Salem set $E\subset \R$ with no 3-APs.
\end{theorem}

From this viewpoint, Shmerkin's result says that when it comes to 3-APs, Fourier dimension is \textit{just as bad} at detecting patterns in a set. On the other hand, the following result of Liang--Pramanik shows that sets of large Fourier dimension cannot avoid having nontrivial solutions to \textit{all} translation invariant linear equations.

\begin{theorem}[Liang--Pramanik {\cite[Theorem 1.2]{LiangPramanik}}] \label{liangpramanikNonAvoidance}
    Any set $E \subset \R$ with $\dim_{\mathbb{F}}(E) > \frac{2}{3}$ must contain a nontrivial solution to an equation of the form 
    \[
    (a+b) z = ax + by
    \]
    for some $a,b\in \mathbb{N}$. More generally, any set $E \subset \R$ with $\dim_{\mathbb{F}} (E) > \frac{2}{v+1}$ must contain a nontrivial solution to an equation of the form 
    \[
    (m_1 + \cdots + m_v) x_0 = m_1x_1 + \cdots + m_v x_v
    \]
    for some $m_i \in \mathbb{N}$.
\end{theorem}

\begin{remark}
In comparison to M\'ath\'e's example in Theorem \ref{matheFulldim} it follows that any rationally independent set must have null Fourier dimension.
\end{remark}

\noindent We discuss Theorem \ref{liangpramanikNonAvoidance} in Section \ref{utilizingFourierDecay}, but first develop intuition regarding Fourier dimension.

\subsection{Randomness} \label{sec:randomness}

A measure $\mu$ satisfying the Fourier decay condition
\[
|\hat{\mu}(\xi)|\lesssim \langle \xi\rangle^{-s/2},
\]
is said to have $s$-dimensional \textit{randomness}. The goal of this section is to motivate this terminology.

\begin{exercise} \label{SphericalDecayExercise}
    The sphere is sometimes referred to as a pseudo-random object. Motivate this by showing $\mathbb{S}^{n-1} \subset \R^n$ is Salem. In particular, show that the spherical measure in $\R^n$, $\sigma$, satisfies 
    \[
    |\hat{\sigma}(\xi)| \lesssim |\xi|^{-(n-1)/2}
    \]
    by a standard stationary phase argument (cf. \cite[Proposition 5.3]{tao}). Additionally, show by scaling that if $\sigma_s$ is the surface measure on the sphere of radius $s$ centered at the origin, then $\widehat{\sigma_s}(\xi) = s^{n-1}\widehat{\sigma}(s\xi)$; conclude that 
    \begin{equation} \label{eqn:sigma-s-decay}
    |\widehat{\sigma_s}(\xi)|\lesssim s^{(n-1)/2} | \xi|^{- (n-1)/2}.
    \end{equation}
\end{exercise}

We now do a random Cantor construction to obtain a set of any desired Fourier dimension in $[0,1]$. To this end, fix $L_1\leq L_2 \leq \cdots$ and $N_1\leq N_2 \leq \cdots$ such that $L_j \leq N_j$ for all $j\geq 1$. We now begin our Cantor set construction, starting of course with $E_0 = [0,1]$. To construct $E_1$, we break $[0,1]$ into $N_1$ equal intervals of length $\frac{1}{N_1}$. Then, we choose, in expectation, $L_1$-many intervals indexed by $X_1^{(1)} \subset \{0,1,\dots, N_1-1\} := \mathbb{Z}/N_1 \mathbb{Z}$ randomly. In particular, for each $a\in \mathbb{Z}/N_1 \mathbb{Z}$, suppose $a\in X_1^{(1)}$ with probability $\frac{L_1}{N_1}$. Then, we let 
\[
E_1 = \bigcup_{j_1} \{\frac{j_1}{N_1} + [0,\frac{1}{N_1}] : j_1\in X_1^{(1)}\}.
\]
To construct $E_2$, for each $j_1 \in \mathbb{Z}/N_1 \mathbb{Z}$, consider a collection of i.i.d. random subsets $X_{j_1}^{(2)} \subset \mathbb{Z}/N_2\mathbb{Z}$ such that 
\[
\#X_{j_1}^{(2)} = L_2, \quad \mathbb{P}(a \in X_{j_1}^{(2)}) = \frac{L_2}{N_2}, \quad \text{ and }
 \{X_{1}^{(1)}, X_{j_1}^{(2)}\} \text{ are independent.}
 \]
Then, we let 
\[
E_2 = \bigcup_{j_1,j_2} \{\frac{j_1}{N_1} + \frac{j_2}{N_1N_2} + [0,\frac{1}{N_1N_2}] : j_1 \in X_1^{(1)}, j_2 \in X_{j_1}^{(2)}\}.
\]
The construction proceeds iteratively, with the sets $X^{(1)}_1$, $X_{j_1}^{(2)}$, $X_{j_1,j_2}^{(3)},\dots$ being chosen independently. This gives rise to Fourier decay via the following lemma.

\begin{lemma}[Shmerkin {\cite[Theorem 2.1]{Shmerkin2016}}] \label{shmerkin-lemma}
    Assume 
    \[
    \lim_{k\to \infty} \frac{\log N_{k+1}}{ \log N_1\cdots N_k} = 0.
    \]
    Then, letting $E_k$ be as before and $E = \bigcap E_k$, almost surely
    \[
    \dim_{\mathbb{H}}(E) = \dim_{\mathbb{F}}(E) = \liminf_{k\to \infty} \frac{\log (L_1\cdots L_k)}{\log(N_1\cdots N_k)} =: s.
    \]
\end{lemma}

\begin{exercise}
As in Exercise \ref{exercise:h-dimBehrend}, prove $\dim_{\mathbb{H}}(E) =:s$ in the above construction. 
\end{exercise}

\begin{remark}
It is generally not easy to come up with the right strategy for proving Fourier decay of a random measure. Usually, this proceeds by taking a sequence of measures $\mu_k$ supported on $E_k$ such that $(\mu_k - \mu_{k-1})\,\,\hat{}$ is quantitatively small almost surely, and such that $\mu_k$ converges to a measure $\mu$ witnessing the proper Fourier decay. For reference, we refer the reader to \cite[Theorem 2.1]{Shmerkin2016}.
\end{remark}

\begin{exercise} \label{ShmerkinExercise}
Modifying the proof of Theorem~\ref{ctmBehrend}, use a random Cantor set construction to obtain a Salem set $E \subset [0,1]$ of full dimension that avoids 3-term arithmetic progressions. The following outline may be useful. Throughout, identify $[0,1)$ with $\mathbb{R}/\mathbb{Z}$.

\begin{enumerate}
    \item \emph{A translation-stable discrete input.} For $m \in \mathbb{N}$ and $j \in \mathbb{Z}/m\mathbb{Z}$, let $I_{m,j} = [\frac{j}{m}, \frac{j+1}{m})$. Using Theorem~\ref{Behrend}, show that for every $\epsilon > 0$ and all sufficiently large $m$, there exists $X \subset \mathbb{Z}/m\mathbb{Z}$ with $\#X \geq m^{1-\epsilon}$ satisfying the following: for every $\ell \in \mathbb{Z}/m\mathbb{Z}$, if $\{x,y,z\} \subset \bigcup_{j \in X+\ell} I_{m,j}$ is a 3-AP, then all three points lie in a single interval $I_{m,j}$. \emph{Hint:} Take $X = 2B$, where $B \subset \{0, \dots, \lfloor m/5 \rfloor - 1\}$ is 3-AP free, and explain why it suffices to consider $\ell = 0$.

    \item \emph{The random construction.} Let $X^{(k)} \subset \mathbb{Z}/N_k\mathbb{Z}$ be a set of maximal size with the property in Step (1) (with $m = N_k$), and let $L_k = \#X^{(k)}$. Starting from $E_0 = [0,1]$, construct $E_k$ from $E_{k-1}$ as follows: divide each basic interval $I$ of $E_{k-1}$ into $N_k$ equal pieces, and keep those indexed by $X^{(k)} + \ell_I \pmod{N_k}$. Here, the $\ell_I$ are uniform on $\mathbb{Z}/N_k\mathbb{Z}$ and independent across all basic intervals and all stages. Check that exactly $L_k$ pieces are kept (deterministically), and that each individual piece is kept with probability $L_k/N_k$.

    \item \emph{Dimension.} Taking $N_1 = 2$ and $N_k = k$ for $k \geq 2$, use Lemma~\ref{shmerkin-lemma} to show that $E = \bigcap E_k$ almost surely satisfies $\dim_{\mathbb{H}}(E) = \dim_{\mathbb{F}}(E) = 1$.

    \item \emph{No 3-APs.} Show that almost surely $E$ contains no endpoint of a basic interval. Then, given a 3-AP in $E$, consider the smallest basic interval containing all three points, and use Step (1) to reach a contradiction.
\end{enumerate}
For a full proof, see \cite[Theorem B]{Shmerkin2016}. More generally, using Lemma~\ref{lem:partial-avoid} in place of Step 1, try to show that for any finite $\mathcal{G} \subset \mathcal{F}$ there exists a Salem set $E \subset [0,1]$ of full dimension avoiding nontrivial solutions to every $f \in \mathcal{G}$; see \cite[Theorem 1.5]{LiangPramanik}.
\end{exercise}

\subsection{Utilizing Fourier Decay} \label{utilizingFourierDecay}

We now develop some tools from harmonic analysis which utilize the Fourier decay of a measure with $s$-dimensional randomness to detect point configurations in a given set. For the sake of illustration, let us prove the first part of Theorem \ref{liangpramanikNonAvoidance}. That is: any set $E\subset \R$ with $\dim_{\mathbb{F}}(E) > \frac{2}{3}$ must contain a nontrivial solution to an equation of the form 
\[
(a + b) z = ax + by
\]
for some $a,b\in \mathbb{N}$. Notice that since $a,b$ are positive integers, it suffices to show that for some $t\in (0,1)$, $E$ contains a nontrivial solution to an equation of the form 
\begin{equation} \label{tequation}
z = t x + (1-t) y.
\end{equation}
To this end, we construct an integral to detect such configurations. 

If $E$ had positive measure, it would suffice to show that for some $f\geq 0$ with $\mathrm{supp}(f) \subset E$, 
\[
\iint f(x) f(y) f(tx + (1-t)y) \,\mathrm dx \,\mathrm dy >0.
\]
However, since $E$ is a fractal set, we hope to consider a similar trilinear form replacing $f$ with a probability measure $\mu$. A natural choice of probability measure is one witnessing the appropriate Fourier decay of $E$. Hence, fix a probability measure $\mu$ with $\mathrm{supp}(\mu) \subset E$ satisfying 
\begin{equation} \label{lem-fourier-decay}
|\hat{\mu}(\xi)| \lesssim \langle \xi\rangle^{-(\frac{1}{3} + \epsilon_0)}.
\end{equation}
where $\dim_{\mathbb{F}}(E) > \frac{2}{3} + 2\epsilon_0$. We now need make sense of a singular integral of the form 
\[
\iint \mu(x) \mu(y) \mu(tx + (1-t)y)\,\mathrm dx \,\mathrm dy.
\]
To do so, we will need to convolve $\mu$ with an approximate identity $\phi_\epsilon$ (see \cite[Section 1.2.4]{Grafakos2014}); let $\mu_\epsilon = \mu \ast \phi_\epsilon$. We then let
\[
\Lambda_t(\mu) = \lim_{\epsilon \to 0} \iint \mu_\epsilon(x) \mu_\epsilon(y) \mu_\epsilon(tx + (1-t)y) \,\mathrm dx \,\mathrm dy.
\]
This immediately begs two questions.

\begin{question}\label{lambdatquestions}
Why is $\Lambda_t(\mu)$ well-defined, and what does it represent?
\end{question}

In regards to the former part of Question \ref{lambdatquestions}, we claim that $\Lambda_t(\mu)$ can be thought of as a measure due to the Fourier decay of $\mu$. To see this, we can make sense of $\Lambda_t(\mu)$ as a bounded linear functional and apply the Riesz representation theorem. To this end, let 
\[
\langle \Lambda_t(\mu), f\rangle = \lim_{\epsilon \to 0} \iint \mu_\epsilon(x) \mu_\epsilon(y) \mu_\epsilon(tx + (1-t)y) f(x,y)\,\mathrm dx\,\mathrm dy
\]
where $f\in C([0,1]^2).$

\begin{lemma} \label{Fourierform+cty}
    Given $\mu$ satisfying \eqref{lem-fourier-decay}, 
    \[
    \langle \Lambda_t(\mu), 1\rangle = \int_{\R} \hat{\mu}(-t\xi) \hat{\mu}(\xi) \hat{\mu}(-(1-t)\xi) \,\mathrm d\xi < \infty.
    \]
    Moreover, $\Lambda_t(\mu)$ is continuous in $t$ for $t \in (0,1)$.
\end{lemma}

\begin{proof}
    By the Fourier inversion formula and the definition of $\Lambda_t(\mu)$, we have 
    \begin{align*}
    \langle \Lambda_t(\mu), 1\rangle &= \lim_{\epsilon \to 0} \iint \mu_\epsilon(x) \mu_\epsilon(y) \mu_\epsilon(tx + (1-t)y) \,\mathrm dx\,\mathrm dy \\
    &= \lim_{\epsilon \to 0}\int_{\R^2} \int_{\R^3} \hat{\mu_\epsilon}(\zeta) e^{2\pi i x\zeta} \hat{\mu_\epsilon}(\eta) e^{2\pi i y \eta} \hat{\mu_\epsilon}(\xi) e^{{2\pi} i (tx + (1-t) y) \xi} \,\mathrm d\zeta \,\mathrm d\eta \,\mathrm d\xi\,\mathrm dx \,\mathrm dy.
    \end{align*}
    Integrating in $x,y$, we see that the integral vanishes unless $\zeta = -t\xi$ and $\eta = -(1-t)\xi$. Therefore,
    \[
    \langle \Lambda_t(\mu), 1\rangle = \lim_{\epsilon \to 0} \int \hat{\mu_\epsilon}(-t\xi) \hat{\mu_\epsilon}(\xi) \hat{\mu_\epsilon}(-(1-t) \xi) \,\mathrm d\xi.
    \]
    By the dominated convergence theorem, using that $\hat{\mu_\epsilon} = \hat{\mu} \hat{\phi_\epsilon}$ (and $\hat{\phi_\epsilon}$ weakly converges to a Dirac delta), we see that the limit exists and
    \[
    \langle \Lambda_t(\mu), 1\rangle = \int \hat{\mu}(-t\xi) \hat{\mu}(\xi) \hat{\mu}(-(1-t) \xi) \,\mathrm d\xi< \infty.
    \]
    Note that the last inequality holds as $\hat{\mu} \in L^3(\R)$ (see Exercise \ref{hatmuLp}). Finally, one can also glean from the above equation that $\Lambda_t(\mu)$ is continuous in $t$ as claimed.
\end{proof}

\begin{exercise}
    Generalizing the above proof to general $f\in C([0,1]^2),$ show that 
    \[
    |\langle \Lambda_t(\mu), f\rangle| \leq C(t) \lVert f\rVert_\infty.
    \]
\end{exercise}

Hence, by the Riesz representation theorem, it follows that $\Lambda_t(\mu)$ is a measure. Moreover, one can see that $\Lambda_t(\mu)$ is a nonnegative measure given $\mu$ is a probability measure. We can now address the latter part of Question \ref{lambdatquestions}: $\Lambda_t$ counts the number of configurations $x,y,z$ satisfying \eqref{tequation}. In particular, we claim that the support of $\Lambda_t(\mu)$ is contained in the set of solutions to \eqref{tequation} in $E$.

\begin{exercise} \label{exerciseLambdanontrivial}
    Show that
    \[
    \mathrm{supp} (\Lambda_t(\mu)) \subset \{(x,y) \in E \times E : tx + (1-t) y\in E\}.
    \]
    Moreover, show $\Lambda_t(\{(x,x) : x\in E\}) = 0$. Deduce that if $\lVert \Lambda_t(\mu)\rVert >0$, then $E$ contains nontrivial solutions of $tx + (1-t) y = z.$
\end{exercise}

We are now ready to prove the first part of Theorem \ref{liangpramanikNonAvoidance}.

\begin{proof}[Proof of Theorem \ref{liangpramanikNonAvoidance}, first part] Let $\mu$ be as before, and suppose $E$ has no nontrivial solutions to equations of the form 
\[
z = tx + (1-t) y \quad \text{ for any } t\in (0,1) \cap \mathbb{Q}.
\]
Then, by Exercise \ref{exerciseLambdanontrivial}, $\Lambda_t(\mu) \equiv 0$ for all such $t$. By the continuity of $\Lambda_t(\mu)$ in $t$, it follows that $\Lambda_t(\mu) \equiv 0$ for all $t\in (0,1).$ We will use this to reach a contradiction with the fact that $\mu$ is a probability measure. In particular, we have
\begin{align*}
    1 &= \int_{\R^3} \mu(x) \mu(y) \mu(z) \,\mathrm dx\,\mathrm dy\,\mathrm dz \\
    &:= \lim_{\epsilon \to 0} \int_{\R^3} \mu_\epsilon(x) \mu_\epsilon(y) \mu_\epsilon(z) \,\mathrm dx\,\mathrm dy\,\mathrm dz \\
    &= \lim_{\epsilon \to 0} 3! \int_{x < z < y}\mu_\epsilon(x) \mu_\epsilon(y) \mu_\epsilon(z) \,\mathrm dx\,\mathrm dy\,\mathrm dz.
\intertext{Proceeding with the change of variables $z = tx + (1-t) y$, we see
}
    &= \lim_{\epsilon \to 0} 3!\int_{x<y} \int_{t}  \mu_\epsilon(x) \mu_\epsilon(y) \mu_\epsilon(tx + (1-t) y)  |x-y|\,\mathrm dx\,\mathrm dy\,\mathrm dt \\
    &= 3 \int_{[0,1]} \langle \Lambda_t(\mu), |x-y|\rangle \,\mathrm dt = 0
\end{align*}
which gives a contradiction.
\end{proof}

\begin{exercise}
    Sketch the proof that if $\dim_{\mathbb{F}}(E) > \frac{2}{v+1}$, then $E$ contains nontrivial solutions to equations of the form 
    \[
    (m_1 + \cdots + m_v) x_0 = m_1x_1 + \cdots + m_v x_v.
    \]
    In particular, write down the analogue of $\Lambda_t(\mu)$ in this setting and show it's well-defined using the Fourier decay of an appropriately chosen $\mu$. For careful details, see \cite[Proposition 2.1]{LiangPramanik}.
\end{exercise}

\section{Abundance} \label{sec:abundance}

In the following two sections, we aim to give an introduction to a number of other problems in geometric measure theory from the background developed thus far. This section will be on the topic of \textit{abundance}. Throughout, we let $|\cdot|$ denote Lebesgue measure.

Observe that as a corollary of the proof of Theorem \ref{liangpramanikNonAvoidance}, we obtain the following.

\begin{corollary} \label{LiangPramanikPosMeasureCor}
    Given any set $E\subset \R$ with $\dim_{\mathbb{F}}(E)>\frac{2}{3}$, let 
    \[
    T := \{t \in [0,1] : \text{ there exists $x,y\in E$ distinct such that } tx + (1-t)y \in E\}.
    \]
    Then, $|T|>0$.
\end{corollary}

This is an \textit{abundance} result---saying not only does $E$ contain \textit{a} nontrivial solution to an equation of the form \eqref{tequation} if it is sufficiently large, but in fact contains \textit{many}. Abundance problems are, for lack of a better word, abundant in this area of fractal geometry, and in this section we give an introduction to some of them.

\subsection{Falconer-type Abundance} Perhaps one of the most well-known abundance problems is Falconer's distance conjecture. Given $X\subset \R^n$, define the \textit{distance set of} $X$ as 
\[
\Delta(X) := \{|x-y|: x,y\in X\} \subset \R.
\]
In 1985, Falconer asked: How large does $\dim_{\mathbb{H}}(X)$ need to be such that $\Delta(X)$ has positive Lebesgue measure? In particular, he proved the following.

\begin{theorem}[Falconer {\cite{Falconer85}}] \label{thm:falconerfail<n/2}
    Let $n\geq 2$. Then, for all $s< n/2$, there exists a compact set $E\subset \R^n$ with $\dim_{\mathbb{H}}(E) = s$ and $|\Delta(E)| = 0.$
\end{theorem}

\begin{proof}
This result uses a discrete grid-like example. To this end, consider the discrete set $E_q = \mathbb{Z}^n \cap [0,q]^n$. Then, $\# E_q \sim q^n$. Moreover, 
\[
\# \Delta(E_q) \leq \#\{x_1^2 + \cdots + x_n^2 : x_j \in [0,q]\} \lesssim_n q^2 \sim (\# E_q)^{2/n}.
\]
We now put a small ball around each lattice point, and take the union of these balls. To this end, let $q_1 = 2$, $q_{i+1} = q_i^{i+1}$, and let $E_i$ be the $q_i^{-\frac{n}{s}}$ neighborhood of $\frac{1}{q_i}\cdot E_{q_i}$. Finally, let $E = \bigcap E_{i}$. By \cite[Chapter 8]{FalconerBook}, $\dim_{\mathbb{H}}(E) = s$, while 
\[
|\Delta(E_i)| \lesssim q_i^{-\frac{n}{s}} \#\Delta(E_{q_i}) \lesssim q_i^{2-\frac{n}{s}}.
\]
Thus, if $s< \frac{n}{2}$, sending $i$ to infinity implies that $|\Delta(E)| = 0$, proving the claim.
\end{proof}

\begin{remark}
    The discrete example is closely related to the Erd\H{o}s distance problem---a discrete analogue of the Falconer distance problem which was posed by Erd\H{o}s in 1945. This problem asks: Given $E\subset \R^n$ a finite set, then for every $\epsilon>0$, there exists $C_\epsilon>0$ such that
    \[
    \# \Delta (E) \geq C_\epsilon (\# E)^{\frac{2}{n} - \epsilon}.
    \]
    This discrete example shows that the exponent here is sharp. The problem was solved by Guth and Katz in $\R^2$ \cite{GuthKatz}   and very recently by Tidor, Yu, and Zakharov in $\R^3$ \cite{distancesR3}. In higher dimensions, the problem is still open.
\end{remark}

Falconer also proved the following sufficiency condition in the same paper.

\begin{theorem}[Falconer {\cite{Falconer85}}] \label{thm:falconern+1/2}
    Let $n\geq 2$. Given $E\subset \R^n$ compact with $\dim_{\mathbb{H}}(E) > \frac{n+1}{2}$, $|\Delta(E)|>0$.
\end{theorem}

\begin{proof}
We follow Falconer's original proof using an $L^\infty$ method outlined nicely in \cite{IosevichWhatIsFalconer}.

By \eqref{EnergyFrostman}, let $\mu$ be a probability measure with $\mathrm{supp}(\mu) \subset E$ such that $I_{(n+1)/2}(\mu)< \infty.$ Moreover, letting $\phi_\epsilon$ be an approximate identity and $\sigma_s^\epsilon := \sigma_s \ast \phi_\epsilon$, let 
\[
\de(\mu)^\epsilon(s) = \iint \sigma_s^\epsilon(x-y)\,\mathrm d\mu(x)\,\mathrm d\mu(y).
\]
We now want to show that the Lebesgue measure of $\Delta(E)$ is positive. Hence, let $\Delta(E)$ be covered by an arbitrary collection of $\epsilon_i$ intervals $\{(s_i-\epsilon_i, s_i + \epsilon_i)\}.$ Then, we have 
\[
    1 = \mu \times \mu(E \times E)
    \leq \sum_i \mu \times \mu \{(x,y) : s_i - \epsilon_i \leq |x-y|\leq s_i + \epsilon_i\} 
\]
Plugging in the definition of $\de(\mu)^{\epsilon_i}$, we see
\[
1\lesssim \sum_i \epsilon_i \de(\mu)^{\epsilon_i} (s_i) \\
    := \sum_i \epsilon_i \iint \sigma_{s_i}^{\epsilon_i} (x-y) \,\mathrm d\mu(x) \,\mathrm d\mu(y) = \sum_i \epsilon_i \int |\hat{\mu}(\xi)|^2 \hat{\sigma_{s_i}}^{\epsilon_i}(\xi) \,\mathrm d\xi.
\]

\begin{exercise} Using \eqref{eqn:sigma-s-decay}, show
\[
|\hat{\sigma_{s}}^{\epsilon}(\xi)| \lesssim s^{(n-1)/2}|\xi|^{-(n-1)/2}
\]
uniformly in $\epsilon$.
\end{exercise}

Note that by the compactness, one can take $s\leq R$ where $E\subset B^n(0,R)$. Thus, plugging this exercise back into the above and using the compactness of $E$, we have 
\[
    1 \lesssim_R \sum_i \epsilon_i \int |\xi|^{-(n-1)/2} |\hat{\mu}(\xi)|^2 \,\mathrm d\xi = I_{(n+1)/2}(\mu) \sum_i \epsilon_i \lesssim \sum_i \epsilon_i. 
\]
Since the covering was arbitrary, it follows that $|\Delta(E)|>0.$
\end{proof}

Given Theorems \ref{thm:falconerfail<n/2} and \ref{thm:falconern+1/2}, Falconer conjectured the following open problem.

\begin{conjecture}[Falconer] \label{conj:falconer}
    Let $n\geq 2$. If $E\subset \R^n$ is compact with $\dim_{\mathbb{H}}(E) > \frac{n}{2}$, then $|\Delta(E)| >0$.
\end{conjecture}

\begin{remark}
    Note that Conjecture \ref{conj:falconer} fails if $n = 1$; in fact, there exists a compact set $E\subset [0,1]$ with $\dim_{\mathbb{H}}(E) = 1$ such that $|\Delta(E)| = 0$ (see \cite[Section 4.12]{Mattila95}).
\end{remark}

We refer the reader to the best results in the plane due to Guth--Iosevich--Ou--Wang \cite{GuthIosevichOuWang}, and in higher dimensions due to Du--Ou--Ren--Zhang \cite{DuOuRenZhang}.

It is also worth discussing \textit{multipoint} configuration problems. While the Falconer distance problem regards the abundance of a given pattern ($|x-y|$) between \textit{pairs} of points, we can also naturally study configurations with more geometric complexity.

\begin{theorem}[Greenleaf--Iosevich {\cite{GreenleafIosevich}}] \label{GreenleafIosevich}
    Any $E\subset \R^2$ with $\dim_{\mathbb{H}} E>\frac{7}{4}$ spans a positive measure set of triangles. That is, 
    \[
    T(E) := \{(|x-y|,|y-z|,|x-z|) : x,y,z\in E\} \subset \R^3
    \]
    satisfies $|T(E)|>0$.
\end{theorem}

Greenleaf and Iosevich's approach uses an $L^\infty$ method similar to Falconer's original approach with a corresponding trilinear convolution operator. This (Fourier) integral operator approach to more general multipoint configuration problems was generalized by Greenleaf, Iosevich, and Taylor in \cite{GIT1,GIT3,GIT2}.

\begin{exercise}
    Using \cite[(1-2)]{GreenleafIosevich}, which states that for $\mu$ $s$-Frostman for some $s>\frac{7}{4}$,
    \[
    (\mu \times \mu \times \mu)\{(x_1,x_2,x_3) : s_{ij}^\ell - \epsilon_\ell \leq |x^i - x^j | \leq s_{ij}^\ell + \epsilon_\ell \,\forall i,j\} \lesssim \epsilon^3,
    \]
    as well as \cite[(2-3)]{GreenleafIosevich} which gives (for $s = (s_{12},s_{13},s_{23})$)
    \[
    \Lambda_s^{\epsilon}(\mu) = \iiint \sigma_{s_{12}}^{\epsilon} (x^1 - x^2) \sigma_{{s_{13}}}^{\epsilon}(x^1 - x^3) \sigma_{s_{23}}^\epsilon (x^2 - x^3) \,\mathrm d\mu(x^1)\,\mathrm d\mu(x^2) \,\mathrm d\mu(x^3) \lesssim 1, 
    \]
    derive Theorem \ref{GreenleafIosevich} by an $L^\infty$ method. For proofs of both inequalities, we refer the reader to \cite{GreenleafIosevich}.
\end{exercise}

\subsection{Mattila--Sj\"olin-type Abundance} Both Theorems \ref{thm:falconern+1/2} and \ref{GreenleafIosevich} show a positive \textit{measure} set of configurations, but we can also ask when our set of configurations has additional \textit{structure}. For instance, Steinhaus showed that given a set $E\subset \R^n$ with positive measure, $E - E := \{x-y : x,y\in E\}$ contains a ball around the origin. In particular, the difference set $E - E$ contains non-empty interior for $E$ sufficiently large. It is thus reasonable to ask if a similar non-empty interior phenomena happens for other point configurations.

A key result to this end is due to Mattila--Sj\"olin who improved upon Falconer's original distance set result (Theorem \ref{thm:falconern+1/2}).

\begin{theorem}[Mattila--Sj\"olin \cite{MattilaSjolin}] \label{MattilaSjolin}
    Let $n\geq 2$. Given $E\subset \R^n$ compact with $\dim_{\mathbb{H}}(E) > \frac{n+1}{2}$, $\Delta(E)$ has non-empty interior.
\end{theorem}

The above non-empty interior result follows from the following continuity property of $\de(\mu)$.

\begin{theorem}[Mattila--Sj\"olin \cite{MattilaSjolin}] \label{MattilaSjolinContinuity}
    Let $n\geq 2$. If $I_{(n+1)/2}(\mu) < \infty$, then $\de(\mu) = g\,\mathrm ds$ where 
    \[
    g(s) = \begin{cases}
        \int_{\R^n} \widehat{\sigma_s}(\xi) |\hat{\mu}(\xi)|^2 \,\mathrm d\xi & s>0 \\
        0 & s \leq 0.
    \end{cases}
    \]
    The function $g$ is continuous on $\R$ and satisfies $|g(s)|\lesssim s^{(n-1)/2} I_{(n+1)/2}(\mu)$.
\end{theorem}

\begin{proof}[Proof of Theorem \ref{MattilaSjolin}]
    Taking $\mu$ a probability measure with $\mathrm{supp}(\mu) \subset E$ and $I_{\frac{n+1}{2}}(\mu)< \infty$, notice 
    \[
    \int_{\R} g(s) \,\mathrm ds = \de(\mu)(\R) = \mu(E)^2=1.
    \]
    So, $g(s_0) >0$ for some $s_0$. In particular, by continuity, $g>0$ on an open neighborhood $U\ni s_0$. Thus, $U \subset \mathrm{supp}(\de(\mu)) \subset \Delta(E).$
\end{proof}

\begin{proof}[Proof of Theorem \ref{MattilaSjolinContinuity}]
    We begin by presenting an alternative proof of Theorem \ref{thm:falconern+1/2} due to Mattila in 1987 \cite{Mattila_1987}. The proof also utilizes Exercise \ref{SphericalDecayExercise}---the Fourier decay of spherical measures. Moreover, we will need to recall the energy-theoretic version of Frostman's lemma \eqref{EnergyFrostman}: Given $\dim_{\mathbb{H}}(E)> \frac{n+1}{2}$, there exists a probability measure $\mu$ with $\mathrm{supp} (\mu ) \subset E$ satisfying 
\[
I_s(\mu):= \int_{\R^n} |\hat{\mu}(\xi)|^2 |\xi|^{s-n} \,\mathrm d\xi < \infty.
\]
Finally, we will let $\de(\mu)$ denote the pushforward of $\mu\times \mu$ under the distance function, i.e.
\[
\de(\mu)(A) := (\mu \times \mu) \{(x,y) : |x-y| \in A\}.
\]
Notice that since $\mu$ has finite energy, $(\mu \times \mu)(\{x = y\}) = 0$, and in particular $\de(\mu)(0) = 0$. 

\begin{exercise} \label{DensityDeMu}
    Given $f\in C_c^\infty(\R^n)$ real valued, then for $s>0$, 
    \[
    \de(f)(s) = \int_{\R^n} (f\ast \sigma_s)(x) f(x) \,\mathrm dx = \int_{\R^n} \widehat{\sigma_s}(\xi) |\hat{f}(\xi)|^2 \,\mathrm d\xi
    \]
    Note that the last equality holds by Plancherel.
\end{exercise}

Using Exercise \ref{DensityDeMu} and \eqref{eqn:sigma-s-decay}, we see that for $f\in C_c^\infty(\R^n)$.
\begin{align*}
|\de(f)(s)| &\leq \int_{\R^n}  |\widehat{\sigma_s}(\xi)| |\hat{f}(\xi)|^2 \,\mathrm d\xi \\
&\lesssim s^{(n-1)/2} \int_{\R^n} |\xi|^{-(n-1)/2} |\hat{f}(\xi)|^2 \,\mathrm d\xi \\
&\lesssim s^{(n-1)/2} I_{(n+1)/2} (f).
\end{align*}
Notice that given $E\subset \R^n$ is compact, there exists some $R>0$ such that $\de(\mu)(s) = 0$ for all $s\geq R$. Thus, 
\begin{equation} \label{eqn:keycty}
|\de(f)(s)| \lesssim R^{(n-1)/2} I_{(n+1)/2} (f) 
\end{equation}
for all $s>0$. 

From here, we can deduce Falconer's bound (Theorem \ref{thm:falconern+1/2}). Using an approximate identity $\phi_\epsilon$, it follows that $\mu_\epsilon := \mu \ast \phi_\epsilon \in C_c^\infty(\R^n)$ converges to $\mu$ weakly on a compact set. Thus, 
\[
|\de(\mu)(s)|  := \lim_{\epsilon \to 0} |\de(\mu_\epsilon)(s)| \lesssim \lim_{\epsilon \to 0} R^{(n-1)/2} I_{(n+1)/2} (\mu_\epsilon) \leq  R^{(n-1)/2} I_{(n+1)/2} (\mu) := M.
\]
Let us justify that this pointwise limit exists, and that it is indeed the density of $\de(\mu)$. Take $\phi \in C_c^\infty(\R^n)$ with $\phi\geq 0$ and $\int \phi = 1$, and $\phi_\epsilon(x) = \epsilon^{-n}\phi(x/\epsilon)$. Then $\widehat{\mu_\epsilon}(\xi) = \hat{\phi}(\epsilon \xi)\hat{\mu}(\xi)$, where $|\hat{\phi}| \leq \int \phi = 1$ and $\hat{\phi}(\epsilon\xi) \to \hat{\phi}(0) = 1$ as $\epsilon \to 0$. By Exercise \ref{DensityDeMu}, for each fixed $s>0$,
\[
\de(\mu_\epsilon)(s) = \int_{\R^n} \widehat{\sigma_s}(\xi)\, |\hat{\phi}(\epsilon\xi)|^2\, |\hat{\mu}(\xi)|^2 \,\mathrm d\xi.
\]
The integrand converges pointwise to $\widehat{\sigma_s}(\xi)|\hat{\mu}(\xi)|^2$, and by \eqref{eqn:sigma-s-decay} it is dominated, uniformly in $\epsilon$, by
\[
|\widehat{\sigma_s}(\xi)||\hat{\mu}(\xi)|^2 \lesssim s^{(n-1)/2} |\xi|^{-(n-1)/2}|\hat{\mu}(\xi)|^2,
\]
which is integrable since its integral is $s^{(n-1)/2} I_{(n+1)/2}(\mu) < \infty$. By dominated convergence,
\[
\lim_{\epsilon \to 0} \de(\mu_\epsilon)(s) = \int_{\R^n} \widehat{\sigma_s}(\xi) |\hat{\mu}(\xi)|^2 \,\mathrm d\xi = g(s),
\]
so the limit exists for every $s > 0$. The same domination shows $I_{(n+1)/2}(\mu_\epsilon) \leq I_{(n+1)/2}(\mu)$, giving the last inequality above. To see that $g$ is the density of $\de(\mu)$, let $h \in C_c((0,\infty))$. Since $\mu_\epsilon \times \mu_\epsilon \to \mu\times\mu$ weakly, with supports in a fixed compact set, and $(x,y) \mapsto h(|x-y|)$ is bounded and continuous,
\[
\int h(s) \de(\mu_\epsilon)(s) \,\mathrm ds = \iint h(|x-y|)\, \mu_\epsilon(x)\mu_\epsilon(y)\,\mathrm dx\,\mathrm dy \longrightarrow \iint h(|x-y|) \,\mathrm d\mu(x)\,\mathrm d\mu(y) = \int h \,\mathrm d\de(\mu).
\]
On the other hand, $|\de(\mu_\epsilon)(s)| \lesssim s^{(n-1)/2} I_{(n+1)/2}(\mu)$ is bounded on $\mathrm{supp}(h)$ uniformly in $\epsilon$, so dominated convergence (now in $s$) gives $\int h(s) \de(\mu_\epsilon)(s)\,\mathrm ds \to \int h(s) g(s) \,\mathrm ds$. Hence $\de(\mu) = g\,\mathrm ds$ on $(0,\infty)$, and since $\de(\mu)(\{0\}) = 0$, on all of $[0,\infty)$.

As a consequence, it follows that $\de(\mu)$ is absolutely continuous with respect to Lebesgue measure with density at most $M$. In particular, 
\[
1 = \mu(E)^2 = (\mu \times \mu)(E \times E) \leq \de(\mu) (\Delta(E)) \leq M |\Delta(E)|.
\]
Thus, $|\Delta(E)| >0$.

Picking up from \eqref{eqn:keycty}, we have more or less already proven Theorem \ref{MattilaSjolinContinuity} by Exercise \ref{DensityDeMu}. In particular, one needs to check the continuity of $g$. Continuity at 0 follows as (we have shown),
\[
|g(s)| \lesssim s^{(n-1)/2} I_{(n+1)/2}(\mu).
\]
The continuity of $g$ is outlined in the following exercise, completing the proof of Theorem \ref{MattilaSjolinContinuity}.
\end{proof}

\begin{exercise}
Fix $s_0>0$. Find a function $G\in L^1(\R^n)$, independent of $s$, such that 
\[
|\widehat{\sigma}_s(\xi)||\hat{\mu}(\xi)|^2 \leq G(\xi) \quad \text{ for all $s\in [s_0/2,2s_0]$ and a.e. $\xi$.}
\]
Use this to prove $g$ is continuous at $s_0$.
\end{exercise}

For more recent progress on Mattila--Sj\"olin-type results for distance sets and more general multipoint configuration problems, we refer the reader to \cite{BorgesEtAl2,BorgesEtAl1,PalssonRomeroAcosta1,PalssonRomeroAcosta2} and the references therein.

\subsection{Ultimate Abundance} We can ask for even more structure than non-empty interiors. For instance, we could ask: What conditions on $E\subset \R^n$ guarantee there exists a $\lambda_0(E)>0$ with $[\lambda_0,\infty) \subset \Delta(E).$ We refer to this as \textit{ultimate abundance}, and such conditions tend to require a different notion of size. A standard notion of size to this end is \textit{positive upper density}. Recall that given $E\subset \R^n$, $E$ has positive upper density if 
\[
\overline{\de}(E) := \lim_{R\to \infty} \sup_{x\in \R^n} \frac{|E\cap B(x,R)|}{|B(x,R)|} >0.
\]

\begin{theorem}[Bourgain, Falconer--Marstrand, Furstenberg--Katznelson--Weiss \cite{BourgainSzemeredi-type,FalconerMarstrand,FurstenbergKatznelsonWeiss}] \label{thm:posupperdensity}
    Let $n\geq 2$. Any $E\subset \R^n$ with positive upper density contains all sufficiently large distances. That is, there exists some $\lambda_0 (E) >0$ such that 
    \[
    \Delta(E):= \{|x-y| : x,y\in E\} \supset [\lambda_0,\infty).
    \]
\end{theorem}

For a proof of this result, Steinhaus' theorem (Theorem \ref{steinhaus}), and much more on finite Euclidean configurations
in sets of positive density, we refer the reader to \href{https://web.math.pmf.unizg.hr/~vjekovac/files/Lecture_notes_Simons_Budapest.pdf}{Kova\v{c}'s lectures} on the topic. It is reasonable to ask if we can find such ultimate abundance if we weaken the condition on positive upper density. A result of Rice answers this question strongly in the negative.

\begin{theorem}[Rice {\cite[Theorem 1]{Rice2019}}]
    Let $n\geq 2$ and let $\psi: (0,\infty) \to [0,1]$ satisfy $\lim_{R\to \infty} \psi(R) = 0.$ Then, there exists a set $E\subset \R^n$ and a sequence of scales $R_j\nearrow\infty$ such that 
    \[
    \sup_{x\in \R^n} \frac{|E\cap B(x,R_j)|}{|B(x,R_j)|} \geq \psi(R_j) \quad \text{ for all $j\geq 1$}\]
    such that $\Delta(E) \cap \{R_j : j\geq 1\} = \emptyset$.
\end{theorem}

On the other hand, a flexibility result of Pramanik--Raani gives a construction of a set with null upper density and an abundance of distances.

\begin{theorem}[Pramanik--Raani {\cite[Corollary 1.6]{PramanikRaani}}]
    Let $n\geq 2$ and $s \geq n - \epsilon_n$ for some sufficiently small $\epsilon_n >0$. Then, there exists a set $E\subset \R^n$ with $\dim_{\mathbb{H}}(E) \leq s$ such that
    \begin{itemize}
        \item $\Delta(E\setminus B(0,R)) = [0,\infty)$ for all $R\geq 1$, and 
        \item $[0,\infty) \setminus \Delta(E \cap B(0,R)) \supset [a_R,b_R]$ such that $b_R - a_R \nearrow \infty$ as $R\nearrow\infty.$
    \end{itemize}
\end{theorem}

\begin{remark}
    Note that by Theorem \ref{thm:posupperdensity} and Exercise \ref{dimensionbasics},  it must be the case that $\overline{\de}(E) = 0$ in both of the previous two theorems.
\end{remark}

\noindent In other words, outside of any ball of radius $R\geq 1$, $E$ contains \textit{all} distances, and yet \textit{within} any such balls some interval of distances are missing. 

\subsection{Distances and Fourier Decay}
Finally, it is worth asking how distance sets interact with Fourier decay. To this end, Pramanik and Raani explored the role of \textit{Fourier asymptotics} in the distance set problem. Given a probability measure $\mu$ such that $\mu\perp \mathcal H^n$, notice that 
\[
\mathbb{F}_\mu(T) = \int_{|\xi| \leq T} |\hat{\mu}(\xi)|^2 \,\mathrm d\xi \to \infty \quad \text{ as } T\to \infty.
\]
It turns out that the rate of this blow up contains information about $\Delta(E)$ where $E = \mathrm{supp}(\mu)$.

\begin{theorem}[Pramanik--Raani {\cite[Theorem 1.7]{PramanikRaani}}] \label{PramanikRaaniThm}
    Suppose $\mathbb{F}_\mu(T) \sim T^\epsilon$ for some $\epsilon>0$. Then, $\Delta(E)$ intersects \textit{every} geometric sequence of the form $\{a \tau^j : j\geq 1\}$ where $a>0$ and $\tau \in (0,1).$ Moreover, there exists $J_0 := J_0(\epsilon,\tau)$ such that 
    \[
    \#\{\Delta(E) \cap \{a\tau^j : R \leq j \leq R + J_0\}\} \geq 1 \quad \text{ for all } R\geq 1.
    \]
\end{theorem}

\begin{remark}
Theorem \ref{PramanikRaaniThm} can be thought of as an abundance result in that the distance set must intersect not just every geometric sequence, but every geometric sequence infinitely often. 
\end{remark}

What if instead of Fourier asymptotics, we used Fourier dimension? In the 1987 paper of Mattila \cite{Mattila_1987}, he showed that given a Borel set $E\subset \R^n$ with 
\begin{equation}\label{mattila-FourierDistance}
\dim_{\mathbb{F}} E + \dim_{\mathbb{H}}(E) > n, \quad \text{ then } |\Delta(E)|>0;
\end{equation}
see Exercise \ref{Exercise:MarstrandFourier}. Thus, Salem sets satisfy Falconer's conjecture (Conjecture \ref{conj:falconer}). 

\begin{exercise}[Fraser {\cite[Theorem 7.3]{Fraser2024}}] \label{Fraser-Exercise}
    Let $n\geq 2$ and $E\subset \R^n$ with $\dim_{\mathbb{F}}(E) > n/2$. By Exercise \ref{hatmuLp},  there exists a probability measure $\mu$ with $\mathrm{supp}(\mu) \subset E$ such that $\hat{\mu}\in L^4(\R^n).$ Show that if $\tilde{\mu}(A) := \mu(-A)$, the measure $\nu = \mu \ast \tilde{\mu} \in L^2(\R^n)$ by Plancherel. Using that $\mu \ast \tilde{\mu}$ has support within $E- E$, show that $E-E$ has positive Lebesgue measure. Therefore, 
    \[
    0 < \int_{\R^n} \mathbf{1}_{E-E}(x) \,\mathrm dx = \int_{\mathbb{S}^{n-1}} \int_0^\infty \mathbf{1}_{E-E}(r\theta) r^{n-1} \,\mathrm dr \mathrm d\sigma(\theta).
    \]
    Conclude that there exists a $\theta \in \mathbb{S}^{n-1}$ such that $|\{r: r\theta \in E- E\}|>0$, and thus $|\Delta(E)|>0.$
\end{exercise}

More recent progress has been made utilizing the \textit{Fourier spectrum}, which is a continuously parameterized family of dimensions interpolated between the Fourier and Hausdorff dimension. We elaborate more upon the Fourier spectrum in the following section, but insofar as it relates to abundance via the distance set problem, we refer the reader to \cite[Section 7]{Fraser2024}.

\section{Recent Progress on 3-APs} \label{sec:recentprogress}

In our final section, we discuss some recent progress on finding 3-APs in fractal sets due to de Orellana and Krause \cite{dOKrause}. This section of the notes is based on a lecture by de Orellana at the Institut Mittag-Leffler. To begin, we note some important motivating results.

Firstly, recall that Szemer\'edi's theorem (Theorem \ref{szemeredi}) implies that subsets of the integers with positive upper density contain 3-APs. This result was first proven by Roth \cite{Roth1953} using a method similar to the configuration integrals we saw in Section \ref{lec2}. To this end, we sketch the proof. Fixing a sufficiently large enough $N$ and letting $[N] := \{1,\dots, N\}$, we can focus our attention to $A\cap [N]$ with positive density so that we may employ Fourier analysis over $\mathbb{Z}_N$. Defining the counting form 
\[
\Lambda(A) := \frac{1}{N^2} \sum_{t,x\in [N]} \mathbf{1}_A(x) \mathbf{1}_A(x + t) \mathbf{1}_A(x + 2t)
\]
over $\mathbb{Z}_N$, we then employ Fourier inversion and use that $A$ has $\de$-density to see 
\[
\Lambda(A) = \de^3 + \sum_{\xi \neq 0} \widehat{\mathbf{1}_A}(\xi)\widehat{\mathbf{1}_A}(-2\xi)\widehat{\mathbf{1}_A}(\xi).
\]
The last term will be our error term. If $|\widehat{\mathbf{1}_A}(\xi)|$ is sufficiently small, this implies that the error is small and thus $|\Lambda(A)| \gtrsim \de^3.$ If, on the other hand, $\widehat{\mathbf{1}_A}(\xi)$ is large, this implies that there are points in $A$ with more density, and thus we can ``zoom in'' and iterate this approach.

A continuum analogue of Roth's theorem was given by {\L}aba--Pramanik in 2009 assuming appropriate Frostman and Fourier conditions on a measure $\mu$:

\begin{theorem}[{\L}aba--Pramanik \cite{LabaPramanik}]
    Given $C_1,C_2>0$ and $s>2/3$, there exists an $\epsilon_0(C_1,C_2,s)$ such that the following holds. If $A\subset [0,1]$ is compact and supports a probability measure $\mu$ with 
    \begin{itemize}
        \item $\mu(B(x,r)) \leq C_1 r^{1-\epsilon_0}$ and 
        \item $|\hat{\mu}(\xi)| \leq C_2 |\xi|^{-s/2}$.
    \end{itemize}
    Then, $A$ contains 3-APs.
\end{theorem}

This may be surprising in comparison to Theorem \ref{shmerkin3AP}, which we recall states that there exists a full dimensional Salem set $E\subset \R$ with \textit{no} 3-APs \cite{Shmerkin2016}. To this end, it is important that the $\epsilon_0$ be quantified in the above theorem, as otherwise the theorem \textit{fails}. In particular, $\epsilon_0$ needs to be taken sufficiently small such that $0< \Lambda(\mu)$.

Note that the $2/3$ condition on $s$ is natural for the same reason it is natural in Liang--Pramanik's work. The authors consider a corresponding counting form $\Lambda(f)$ and to make sure it is well-defined and can be thought of as a measure, they needed $s>2/3$. That said, recent work of de Orellana and Krause generalized this to all $s>0$.

\begin{theorem}[de Orellana--Krause {\cite[Theorem 1.10]{dOKrause}}] \label{thm:dOKrause}
    Given $C_1,C_2>0$ and $s>0$, there exists an $\epsilon_0(C_1,C_2,s)$ such that the following holds. If $A\subset [0,1]$ is compact and supports a probability measure $\mu$ with 
    \begin{itemize}
        \item $\mu(B(x,r)) \leq C_1 r^{1-\epsilon_0}$ for all $x\in \R$ and $r>0$, and 
        \item $|\hat{\mu}(\xi)| \leq C_2 |\xi|^{-s/2}$.
    \end{itemize}
    Then, $A$ contains 3-APs.
\end{theorem}

\noindent Their proof makes use of the Fourier spectrum, a parameterized family of dimensions $\dim_{\mathbb{F}}^\theta$ developed by Fraser \cite{Fraser2024} with $\dim_{\mathbb{F}}^0  =\dim_{\mathbb{F}}$, and $\dim_{\mathbb{F}}^1$ equaling \textit{Sobolev dimension}.

\begin{remark}
    Sobolev dimension is a close relative of Hausdorff dimension, see \cite[Section 5.2]{Mattila95}, which we denote $\dim_{\mathbb{S}}.$
\end{remark}

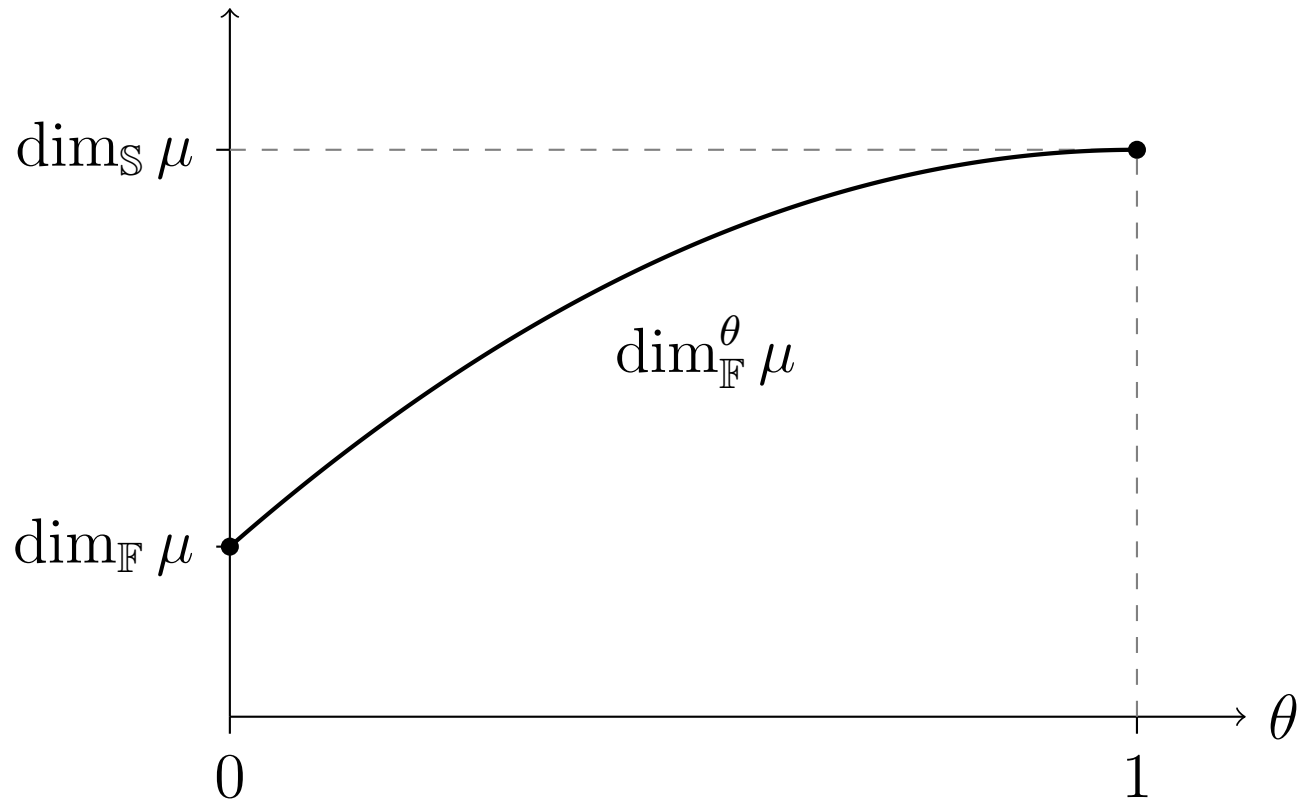
\begin{figure}[ht] 
    \centering
    \begin{tikzpicture}[scale=.75,xscale=8,yscale=5]
  % endpoint values (purely illustrative)
  \pgfmathsetmacro{\a}{0.3}  % = dim_F mu
  \pgfmathsetmacro{\b}{1.0}  % = dim_S mu
  % axes
  \draw[->] (0,0) -- (1.12,0) node[right] {$\theta$};
  \draw[->] (0,0) -- (0,1.25);
  % ticks
  \draw (0,0) -- (0,-0.03) node[below] {$0$};
  \draw (1,0) -- (1,-0.03) node[below] {$1$};
  \draw (0,\a) -- (-0.015,\a) node[left] {$\dim_{\mathbb{F}}\mu$};
  \draw (0,\b) -- (-0.015,\b) node[left] {$\dim_{\mathbb{S}}\mu$};
  % guide lines
  \draw[dashed,gray] (0,\b) -- (1,\b);
  \draw[dashed,gray] (1,0) -- (1,\b);
  % concave, non-decreasing, continuous curve: a + (b-a)(2t - t^2)
  \draw[thick,domain=0:1,samples=100,smooth]
    plot (\x,{\a+(\b-\a)*(2*\x-\x*\x)});
  \node[circle,fill,inner sep=1.2pt] at (0,\a) {};
  \node[circle,fill,inner sep=1.2pt] at (1,\b) {};
  \node[below right] at (0.4,{\a+(\b-\a)*(0.8-0.16)}) {$\dim_{\mathbb{F}}^{\theta}\mu$};
\end{tikzpicture}
\caption{A depiction of $\dim_{\mathbb{F}}^\theta(\mu)$, a non-decreasing, concave, and continuous function in $\theta \in (0,1]$.}
    \label{figure3}
\end{figure}

Analogous to how $s$-energy $I_s(\mu)$ can be used to define Hausdorff dimension \eqref{EnergyFrostman}, the Fourier spectrum is defined by energy integrals $\mathcal J_{s,\theta}$. For $\theta \in (0,1]$ and $\mu$ supported on $A\subset \R$, let 
\[
\mathcal J_{s,\theta}(\mu) := \left(\int_\R |\hat{\mu}(\xi)|^{\frac{2}{\theta}} |\xi|^{\frac{s}{\theta} - 1}\,\mathrm d\xi\right)^\theta,
\]
and let 
\[
\mathcal J_{s,0}(\mu) = \sup_\xi |\hat{\mu}(\xi)|^2 |\xi|^s.
\]
Then, we define the \textit{Fourier spectrum} of $\mu$ at $\theta \in [0,1]$ as 
\[
\dim_{\mathbb{F}}^\theta \mu := \sup \{s\in [0,1] : \mathcal J_{s,\theta}(\mu) < \infty\}.
\]
The Fourier spectrum is a non-decreasing, concave, and continuous in $(0,1]$ (see Figure \ref{figure3}).

Theorem \ref{thm:dOKrause} will follow as a consequence of a more general energy theoretic statement:

\begin{theorem}[de Orellana--Krause {\cite[Theorem 1.11]{dOKrause}}] \label{thm:dOKrauseEnergy}
    Let $A\subset [0,1]$ be compact, and let $A$ support a probability measure $\mu$ with $\mathcal J_{s,\theta}(\mu)< \infty$, $\mathcal J_{u,1-\frac{\theta}{2}}(\mu) < \infty$ for some $\theta \in [0,1]$ and $\frac{s}{2} + u >1 $. Then, given $C_1>0$, there exists an $\epsilon_0(C_1,\mathcal J_{s,\theta}(\mu), \mathcal J_{u,1-\frac{\theta}{2}}(\mu))$  such that if $\mu(B(x,r)) \leq C_1 r^{1-\epsilon_0}$ for all $x \in \R$ and $r>0$, then $A$ contains 3-APs.
\end{theorem}

\begin{remark}
In \cite{dOKrause}, Theorems \ref{thm:dOKrause} and \ref{thm:dOKrauseEnergy} are stated more generally for progressions of the form $\{x, x - \gamma_1 t, x - \gamma_2 t\}$ with $t \neq 0$ and $\gamma_1, \gamma_2$ nonzero rationals, in which case $\epsilon_0$ also depends on $\gamma_1$ and $\gamma_2$; 3-APs correspond to $\gamma_1 = 1$, $\gamma_2 = 2$. Optimizing over $\theta$ in Theorem \ref{thm:dOKrauseEnergy}, they also obtain the following statement, whose numerology matches that of {\L}aba--Pramanik: it suffices that $\mathcal J_{s, 2/3}(\mu) < \infty$ for some $s > \frac{2}{3}$ \cite[Corollary 1.12]{dOKrause}.
\end{remark}

\begin{proof}[Proof of Theorem \ref{thm:dOKrause}]
    Choose $\theta = 0$ and note that since $|\hat{\mu}(\xi)|\leq C_2 |\xi|^{-s/2}$, $\mathcal J_{s,0}(\mu) <C_2^2$. Moreover, we have $\mathcal J_{u,1}(\mu)<\infty$ for any $0 < u < 1-\epsilon_0$ by the following exercise.

    \begin{exercise}
        Given $\mu$ satisfying the Frostman condition $\mu(B(x,r)) \lesssim r^{1-\epsilon_0}$, show $\mathcal J_{u,1}(\mu) < \infty$ for any $0 < u < 1-\epsilon_0$.
    \end{exercise}

    \noindent Noting that one can take $\epsilon_0 < \frac{s}{4}$, we let $u:= 1- \frac{s}{4} < 1-\epsilon_0$. Thus, $\frac{s}{2} + u = 1 + \frac{s}{4}>1$. So, by Theorem \ref{thm:dOKrauseEnergy}, $A$ contains a 3-AP. 
\end{proof}

We spend the rest of the section sketching the proof of Theorem \ref{thm:dOKrauseEnergy}. As before, we will essentially take $\Lambda$ to be the counting form 
\[
\Lambda(\mu) = \int_\R \hat{\mu}(\xi) \hat{\mu}(-2\xi)
 \hat{\mu}(\xi) \,\mathrm d\xi.
 \]
 One needs show that $0 < \Lambda(\mu) < \infty$; the former inequality to see $A$ contains a 3-AP, and the latter to prove $\Lambda(\mu)$ can be thought of as a measure by the Riesz representation theorem. We focus our attention on showing $\Lambda(\mu) < \infty$ as this is what is most different from Lemma \ref{Fourierform+cty} (since $\mu$ does not have large Fourier decay). For the proof that $0 < \Lambda(\mu)$, see \cite[Proposition 4.2]{dOKrause}, and note that this is where we need take $\epsilon_0$ sufficiently small. We stress that finiteness and positivity of the scalar quantity $\Lambda(\mu)$ do not on their own produce a \emph{nontrivial} progression. As with $\Lambda_t(\mu)$ in Section \ref{utilizingFourierDecay}, one must realize $\Lambda(\mu)$ as the total mass of a nonnegative configuration measure that is supported on 3-term progressions in $A$ and assigns zero mass to the diagonal of trivial progressions (those with common difference $0$). In \cite{dOKrause}, this is done in two steps. After the proof of \cite[Proposition 4.2]{dOKrause}, a separation argument shows that $\mathrm{supp}(\mu)^3$ meets the set of progressions, and \cite[Proposition 4.3]{dOKrause} shows, for $\epsilon_0 < \frac12$, that the form assigns no mass to the trivial progressions. We do not reproduce these arguments here.

\begin{proof}[Proof sketch of Theorem \ref{thm:dOKrauseEnergy}]
We proceed by a Littlewood--Paley decomposition. Fix an even, smooth function $\psi \colon \mathbb{R} \to [0,1]$ that is non-increasing in $|\xi|$, with $\psi(\xi) = 1$ for $|\xi| \leq 1$ and $\psi(\xi) = 0$ for $|\xi| \geq 2$. For $\ell \in \mathbb{Z}$, set
\[
\psi_{\leq \ell}(\xi) := \psi(2^{-\ell}\xi), \qquad \varphi_\ell(\xi) := \psi_{\leq \ell}(\xi) - \psi_{\leq \ell-1}(\xi).
\]
Then $0 \leq \varphi_\ell \leq 1$, $\operatorname{supp} \varphi_\ell \subset \{2^{\ell-1} \leq |\xi| \leq 2^{\ell+1}\}$, and for $L < M$ the sum telescopes:
\[
\psi_{\leq L} + \sum_{\ell = L+1}^{M} \varphi_\ell = \psi_{\leq M} \xrightarrow{M \to \infty} 1 \quad \text{pointwise.}
\]
Define $\mu_{\leq \ell}$ and $\mu_\ell$ by their Fourier transforms,
\[
\widehat{\mu}_{\leq \ell} := \psi_{\leq \ell}\, \widehat{\mu}, \qquad \widehat{\mu}_\ell := \varphi_\ell\, \widehat{\mu},
\]
so that $\widehat{\mu} = \widehat{\mu}_{\leq L} + \sum_{\ell > L} \widehat{\mu}_\ell$, with $|\widehat{\mu}_{\leq \ell}|, |\widehat{\mu}_\ell| \leq |\widehat{\mu}|$ pointwise. Note that $\mu_{\leq \ell} = \mu \ast \check{\psi}_{\leq \ell}$ is a smooth function, and $\widehat{\mu}_{\leq \ell}$ is compactly supported. However, $\mu_{\leq \ell}$ need not be a nonnegative function, since the kernel $\check{\psi}_{\leq \ell}$ is not nonnegative. In fact it cannot be: if $\check{\psi} \geq 0$, then for any $0 < |\xi| \leq 1$ we would have $\int (1 - \cos(2\pi x\xi))\,\check{\psi}(x)\,\mathrm dx = \psi(0) - \psi(\xi) = 0$, which forces $\check{\psi} = 0$ almost everywhere. This is harmless for the error term below, which only uses $|\psi_{\leq \ell}|, |\varphi_\ell| \leq 1$, but it means the main term must be handled with some care; see the end of the proof.

For functions $f, g, h$ with compactly supported Fourier transforms, define the trilinear form
\[
\Lambda(f,g,h) := \int_{\mathbb{R}} \widehat{f}(\xi)\, \widehat{g}(-2\xi)\, \widehat{h}(\xi)\, d\xi,
\]
and write $\Lambda(f) := \Lambda(f,f,f)$. We then define
\[
\Lambda(\mu) := \lim_{M \to \infty} \Lambda(\mu_{\leq M}),
\]
whenever the limit exists. Writing $\Lambda_\ell(\mu) := \Lambda(\mu_{\leq \ell}) - \Lambda(\mu_{\leq \ell-1})$, we have the telescoping identity
\[
\Lambda(\mu_{\leq M}) = \Lambda(\mu_{\leq L}) + \sum_{\ell = L+1}^{M} \Lambda_\ell(\mu).
\]
By trilinearity, using $\mu_{\leq \ell} = \mu_{\leq \ell-1} + \mu_\ell$,
\[
\Lambda_\ell(\mu) = \Lambda(\mu_\ell, \mu_{\leq \ell}, \mu_{\leq \ell}) + \Lambda(\mu_{\leq \ell-1}, \mu_\ell, \mu_{\leq \ell}) + \Lambda(\mu_{\leq \ell-1}, \mu_{\leq \ell-1}, \mu_\ell).
\]
Each term contains exactly one factor localized to frequencies $\sim 2^\ell$. Consequently, the integrand of each term is supported in the annulus $A_\ell := \{2^{\ell-2} \leq |\xi| \leq 2^{\ell+1}\}$. Since every factor is pointwise bounded by $|\widehat{\mu}|$, this gives
\begin{equation}\label{eq:LP-error}
|\Lambda_\ell(\mu)| \leq 3 \int_{A_\ell} |\widehat{\mu}(\xi)|^2\, |\widehat{\mu}(-2\xi)|\, d\xi.
\end{equation}
This is the precise form of the heuristic $\Lambda_\ell(\mu) \approx \Lambda(\mu_\ell)$: the error at scale $\ell$ depends only on $\widehat{\mu}$ near frequency $2^\ell$. In particular, if $\sum_{\ell > L} |\Lambda_\ell(\mu)| < \infty$, then $\Lambda(\mu)$ is well defined and
\[
\Lambda(\mu) = \Lambda(\mu_{\leq L}) + \sum_{\ell > L} \Lambda_\ell(\mu).
\]
The first term is the main term, and the sum is the error term. 

We now show that the error term is absolutely summable, and in fact decays exponentially in $L$. To this end, let $\widetilde{A}_\ell := \{2^{\ell-2} \leq |\xi| \leq 2^{\ell+2}\}$, which contains both $A_\ell$ and $-2A_\ell$. The key observation is that finiteness of $\mathcal{J}_{s,\theta}(\mu)$ controls the size of $\widehat{\mu}$ on each such annulus. Indeed, for $s > 0$ and $p \in [2, \infty)$, we have $|\xi|^{\frac{sp}{2}-1} \gtrsim_{s,p} 2^{\ell(\frac{sp}{2}-1)}$ on $\widetilde{A}_\ell$ (regardless of the sign of the exponent), and therefore
\[
\lVert\widehat{\mu}\rVert_{L^p(\widetilde{A}_\ell)}
\lesssim_{s,p} 2^{-\ell(\frac{s}{2} - \frac{1}{p})} \Big( \int_{\mathbb{R}} |\widehat{\mu}(\xi)|^p |\xi|^{\frac{sp}{2}-1} \, d\xi \Big)^{1/p}
= 2^{-\ell(\frac{s}{2} - \frac{1}{p})}\, \mathcal{J}_{s, \frac{2}{p}}(\mu)^{1/2},
\]
where we used that $\frac{s}{2/p} - 1 = \frac{sp}{2} - 1$. When $p = \infty$, the analogous bound $$\lVert\widehat{\mu}\rVert_{L^\infty(\widetilde{A}_\ell)} \lesssim_s 2^{-\ell s/2}\, \mathcal{J}_{s,0}(\mu)^{1/2}$$ follows immediately from the definition of $\mathcal{J}_{s,0}$.

With these annular bounds in hand, we return to \eqref{eq:LP-error}. Applying Hölder's inequality with one factor of $|\widehat{\mu}(\xi)|$ in $L^{2/\theta}$ and the remaining two factors in $L^q$, where $q = \frac{4}{2-\theta}$ so that $\frac{\theta}{2} + \frac{2}{q} = 1$, and noting that $\lVert\widehat{\mu}(-2\,\cdot)\rVert_{L^q(A_\ell)} \leq \lVert\widehat{\mu}\rVert_{L^q(\widetilde{A}_\ell)}$ by a change of variables, we find
\[
|\Lambda_\ell(\mu)| \leq 3\lVert\widehat{\mu}\rVert_{L^{2/\theta}(\widetilde{A}_\ell)}\lVert\widehat{\mu}\rVert_{L^{q}(\widetilde{A}_\ell)}^2.
\]
Since $\frac{2}{q} = 1 - \frac{\theta}{2}$, the annular bounds (applied with exponents $s$ and $u$, respectively) then give
\[
|\Lambda_\ell(\mu)|
\lesssim 2^{-\ell(\frac{s}{2} - \frac{\theta}{2})}\, \mathcal{J}_{s,\theta}(\mu)^{1/2} \cdot 2^{-\ell(u - 1 + \frac{\theta}{2})}\, \mathcal{J}_{u, 1-\frac{\theta}{2}}(\mu)
= 2^{-\ell(\frac{s}{2} + u - 1)}\, \mathcal{J}_{s,\theta}(\mu)^{1/2}\, \mathcal{J}_{u,1-\frac{\theta}{2}}(\mu),
\]
with implicit constant depending only on $s$, $u$, and $\theta$.

Finally, since $\kappa := \frac{s}{2} + u - 1 > 0$ by assumption, summing the resulting geometric series yields
\[
| \Lambda(\mu) - \Lambda(\mu_{\leq L}) | \leq \sum_{\ell > L} |\Lambda_\ell(\mu)| \lesssim \mathcal{J}_{s,\theta}(\mu)^{1/2}\, \mathcal{J}_{u,1-\frac{\theta}{2}}(\mu)\, 2^{-\kappa L}.
\]
In particular, $\Lambda(\mu)$ is well defined, and the main term $\Lambda(\mu_{\leq L})$ approximates it. Thus, to conclude that $\Lambda(\mu) > 0$, it suffices to show that the main term is bounded below by a positive constant independent of $L$ for suitable $L$. This is precisely where the Frostman condition, with $\epsilon_0$ sufficiently small, enters; we now carry this out, following the proof of \cite[Proposition 4.2]{dOKrause}.

Since $\mu_{\leq L}$ may take negative values, however, we cannot directly apply a lower bound for progressions in nonnegative functions to $\Lambda(\mu_{\leq L})$. In the proof of \cite[Proposition 4.2]{dOKrause}, the main term is written as $\lVert \mu_{\leq l_0}\rVert_\infty^3\, \Lambda\big(\mu_{\leq l_0}/\lVert \mu_{\leq l_0}\rVert_\infty\big)$ and \cite[Lemma 4.4]{dOKrause} is applied to the normalized function; this requires the low-frequency kernel to be nonnegative, with mean one (cf. \cite[Section 1.2]{dOKrause}). Since we use the Fourier cutoff $\psi$ for the error analysis, we first replace it by such a nonnegative kernel. Fix an even, real-valued $\chi \in C_c^\infty(\R)$, $\chi \not\equiv 0$, and let
\[
P := \frac{|\check{\chi}|^2}{\lVert \chi \rVert_2^2} \quad \text{and} \quad P_L(x) := 2^L P(2^L x).
\]
Then $P$ is a nonnegative Schwartz function with $\int P = \hat{P}(0) = 1$, and $\hat{P} = (\chi \ast \chi)/\lVert \chi\rVert_2^2$ is smooth, even, compactly supported, and satisfies $|\hat{P}| \leq 1$. Set
\[
f_L := \mu \ast P_L \geq 0 \quad \text{and} \quad \int f_L = 1.
\]

We first check that replacing $\mu_{\leq L}$ by $f_L$ in the main term costs no more than the error term we have already incurred. By definition,
\[
\Lambda(\mu_{\leq L}) - \Lambda(f_L) = \int_{\R} \big(m_\psi(\xi) - m_P(\xi)\big)\, \hat{\mu}(\xi)^2\, \hat{\mu}(-2\xi)\,\mathrm d\xi,
\]
where $m_\psi(\xi) := \psi(2^{-L}\xi)^2\psi(-2^{1-L}\xi)$ and $m_P(\xi) := \hat{P}(2^{-L}\xi)^2\hat{P}(-2^{1-L}\xi)$. Both multipliers are bounded by $1$ and are compactly supported. Moreover, $m_\psi(\xi) = 1$ for $|\xi| \leq 2^{L-1}$, while $|1 - \hat{P}(\eta)| \lesssim \eta^2$ since $\hat{P}$ is smooth and even with $\hat{P}(0) = 1$. Hence
\[
|m_\psi(\xi) - m_P(\xi)| \lesssim \min\{1, (2^{-L}|\xi|)^2\} \quad \text{for all } \xi,
\]
and the difference vanishes for $|\xi| \gtrsim 2^{L}$. Decomposing into the annuli $A_k$ and using the same H\"older and annular bounds as for $\Lambda_\ell(\mu)$, which hold for every $k$, we obtain
\[
|\Lambda(\mu_{\leq L}) - \Lambda(f_L)| \lesssim \mathcal{J}_{s,\theta}(\mu)^{1/2}\mathcal{J}_{u,1-\frac{\theta}{2}}(\mu) \sum_{k \leq L + O(1)} 2^{-2(L-k)}\, 2^{-\kappa k} \lesssim \mathcal{J}_{s,\theta}(\mu)^{1/2}\mathcal{J}_{u,1-\frac{\theta}{2}}(\mu)\, 2^{-\kappa L}.
\]
Here we assumed, as we may after decreasing $s$ and $u$, that $\kappa < 2$. Combined with the error estimate above,
\begin{equation} \label{eq:main-term-comparison}
|\Lambda(\mu) - \Lambda(f_L)| \leq C_0\, 2^{-\kappa L}, \quad \text{where } C_0 \lesssim_{s,u,\theta,\chi} \mathcal{J}_{s,\theta}(\mu)^{1/2}\mathcal{J}_{u,1-\frac{\theta}{2}}(\mu).
\end{equation}

It therefore remains to bound $\Lambda(f_L)$ from below, using that $f_L$ is a bounded, nonnegative function of mass one concentrated near $[0,1]$. For this we pass to physical space: expanding the Fourier transforms as in Lemma \ref{Fourierform+cty}, for Schwartz functions with compactly supported Fourier transforms we have the physical-space formula
\[
\Lambda(f,g,h) = \iint f(x)\, g(x+t)\, h(x+2t)\,\mathrm dx\,\mathrm dt.
\]
We may assume $C_1 \geq 1$, so that $\mu(B(x,r)) \leq C_1 r^{1-\epsilon_0}$ for all $r > 0$ (for $r \geq 1$ this is trivial as $\mu$ is a probability measure). Since $P(z) \lesssim (1 + |z|)^{-2}$, splitting $\int 2^L P(2^L(x-y))\,\mathrm d\mu(y)$ over the shells $|x-y| \sim 2^{k-L}$, $k \geq 0$, gives
\[
\lVert f_L \rVert_\infty \lesssim 2^L \sup_x \sum_{k\geq 0} 2^{-2k}\, \mu(B(x, 2^{k-L})) \lesssim C_1 2^{L\epsilon_0}\sum_{k \geq 0} 2^{-k(1+\epsilon_0)} \leq C_P\, C_1 2^{L\epsilon_0} =: M_L.
\]
This is where the Frostman condition enters. On the other hand, since $\mathrm{supp}(\mu) \subset [0,1]$, letting $I := [-1,2]$ we have $\int_I f_L \geq 1 - \int_{|z| > 2^L} P(z)\,\mathrm dz \geq \frac12$ for $L \geq L_0(\chi)$. This is where the bounded support enters. Thus
\[
G := \frac{f_L \mathbf{1}_I}{M_L} \quad \text{satisfies} \quad 0 \leq G \leq \mathbf{1}_I, \qquad \int_I G \geq \frac{1}{2M_L} = \frac{|I|}{6 M_L}.
\]

Having reduced matters to a bounded nonnegative function on an interval, we can now import a lower bound for progressions from the integer setting. Namely, we use the following continuous form of Roth's theorem with Varnavides averaging (see \cite[Lemma 4.4]{dOKrause}): for every $\delta > 0$ there is $c(\delta) > 0$, which we may take non-decreasing in $\delta$, such that whenever $I$ is an interval and $0 \leq G \leq \mathbf{1}_I$ satisfies $\int G \geq \delta |I|$,
\[
\iint G(x)\, G(x+t)\, G(x+2t)\,\mathrm dx\,\mathrm dt \geq c(\delta)\, |I|^2.
\]
Since $f_L \geq M_L G \geq 0$ pointwise, the physical-space formula gives
\[
\Lambda(f_L) \geq M_L^3 \iint G(x)G(x+t)G(x+2t)\,\mathrm dx\,\mathrm dt \geq 9\, M_L^3\, c\Big(\frac{1}{6M_L}\Big) \geq \frac{1}{24}\, c\Big(\frac{1}{6M_L}\Big),
\]
where in the last step we used $M_L \geq \lVert f_L\rVert_\infty \geq \frac{1}{|I|}\int_I f_L \geq \frac16$.

Finally, we choose the parameters. Let $c_* := \frac{1}{24}\, c\big(\frac{1}{12 C_P C_1}\big) > 0$, which depends only on $C_1$ (and the fixed choice of $\chi$). First choose $L \geq L_0(\chi)$ so large that $C_0 2^{-\kappa L} \leq \frac{c_*}{2}$; this choice depends only on $C_1$, $\mathcal{J}_{s,\theta}(\mu)$, $\mathcal{J}_{u,1-\frac{\theta}{2}}(\mu)$, $s$, $u$, and $\theta$, and not on $\epsilon_0$. Then choose $\epsilon_0 \leq 1/L$, so that $2^{L\epsilon_0} \leq 2$ and hence $M_L \leq 2C_PC_1$. By the lower bound above and the monotonicity of $c(\delta)$, $\Lambda(f_L) \geq c_*$, and so by \eqref{eq:main-term-comparison},
\[
\Lambda(\mu) \geq \Lambda(f_L) - C_0 2^{-\kappa L} \geq \frac{c_*}{2} > 0.
\]
This explains the dependence $\epsilon_0 = \epsilon_0(C_1, \mathcal{J}_{s,\theta}(\mu), \mathcal{J}_{u,1-\frac{\theta}{2}}(\mu))$ in Theorem \ref{thm:dOKrauseEnergy}. Note that the order of the choices matters: the scale $L$ is fixed by the Fourier side, and only afterwards is the Frostman exponent $1-\epsilon_0$ taken close enough to $1$ that $\mu$ looks like a bounded density at scale $2^{-L}$.
\end{proof}

% Removed: now carried out in the proof above.
% \begin{exercise} \label{exercise:dOK}
%     Letting $\mu$ be as in Theorem \ref{thm:dOKrauseEnergy}, show $\Lambda(\mu_{\leq L})$ is bounded by the Frostman condition.
% \end{exercise} 

\begin{exercise}
    As an illustrative case, work through Theorems \ref{thm:dOKrause} and \ref{thm:dOKrauseEnergy} when $\theta = 0.$
\end{exercise}

\begin{exercise}[Fraser {\cite[Theorem 7.3]{Fraser2024}}] \label{Exercise:MarstrandFourier}
    Generalize Exercise \ref{Fraser-Exercise} using H\"older's inequality to show the following. Given $\mu,\nu$ supported on $X,Y \subset \R^n$ respectively, if $\mathcal J_{s,\theta}(\mu)<\infty$ and $\mathcal J_{t,1-\theta}(\nu)<\infty$ for some $s + t > n$, then $|\Delta(X,Y)|>0$ where 
    \[
    \Delta(X,Y) := \{|x-y|: x\in X, y\in Y\}.
    \] As a hint, try to construct a measure supported on $X-Y.$ Conclude that if $\dim_{\mathbb{F}}^{1/2}(\mu) > n/2$, then $|\Delta(X)|>0$. Moreover, taking $\theta = 0$, show one obtains Mattila's result \eqref{mattila-FourierDistance}.
\end{exercise}

\printbibliography

\end{document}